\documentclass{amsart}
\usepackage{amsthm}
\usepackage{amsmath}
\usepackage{amssymb}
\usepackage{mathtools}
\usepackage{url}

\title[Stable intersection of Cantor sets and robust tangency]{
Stable intersection of Cantor sets in higher dimension
 and robust homoclinic tangency of the largest codimension}
\author{Masayuki Asaoka}
\address{Faculty of Science anf Engineering, Doshisha University,
 Kyotanabe 610-0321, Japan.}

\subjclass[2010]{Primary: 37C29, Secondary: 28A80}

\def\RR{\mathbb{R}}
\def\ZZ{\mathbb{Z}}

\def\cU{\mathcal{U}}
\def\cD{{\mathcal{D}}}

\def\vm{{\mathbf{m}}}

\def\ra{{\rightarrow}}

\def\vphi{\varphi}

\def\wh{\widehat}
\def\sh{{\sharp}}

\def\hsp{{\hspace{3mm}}}

\DeclareMathOperator{\Int}{Int}

\DeclareMathOperator{\Mat}{Mat}
\DeclareMathOperator{\Ker}{Ker}

\DeclareMathOperator{\Gr}{Gr}

\DeclareMathOperator{\diam}{diam}

\DeclareMathOperator{\udim}{{\overline{\dim}}}

\newcommand\mc[1]{{\mathcal{#1}}}

\theoremstyle{plain}
\newtheorem{thm}{Theorem}[section]
\newtheorem{prop}[thm]{Proposition}
\newtheorem{lemma}[thm]{Lemma}
\newtheorem{cor}[thm]{Corollary}

\newtheorem*{thmA}{Theorem A}
\newtheorem*{thmB}{Theorem B}

\theoremstyle{definition}

\newtheorem{qtn}[thm]{Question}

\theoremstyle{remark}
\newtheorem{rmk}[thm]{Remark}
\newtheorem{expl}[thm]{Example}

\begin{document}
\maketitle

\begin{abstract}
In this paper, we construct
\begin{enumerate}
 \item  a pair of two regular Cantor sets in higher dimension
 which exhibits $C^1$-stable intersection, and
 \item a hyperbolic basic set
 which exhibits $C^2$-robust homoclinic tangency of the largest codimension
 for any higher dimensional manifold
\end{enumerate}
 using {\it blenders}.
The former implies that an analog of Moreira's theorem on Cantor sets
 in the real line does not hold in higher dimension.
The latter solves a question posed by Barrientos and A.Raibekas.

\end{abstract}

\section{Introduction}
\label{sec:introduction}

\subsection{Blenders}
In \cite{BD96}, Bonatti and D\'iaz introduced a dynamical  mechanism
 called a {\it blender},
 which provides a variety of $C^1$-robust phenomena; 
 persistently non-hyperbolic and  topologically transitive
 diffeomorphisms (\cite{BD96}), 
 robust heterodimensional cycles (\cite{BD08,BDK12,BKR14}), 
 robust homoclinic tangencies (\cite{BD12}),
 robust existence of non-hyperbolic ergodic measures (\cite{BBD16}),
 and
 large robustly topologically mixing sets in symplectic setting (\cite{NP12}).
Berger \cite{Be16} applied blenders in the jet space
 to gave an open set of diffeomorphisms in which
 typical one admits infinitely many periodic attaractors.
Barrientos and Raibekas \cite{BR18} also applied blenders
 in the Grassmannian bundle to
 provide $C^2$-robust tangency of large codimension.

In \cite{BD12}, Bonatti and D\'iaz gave a formal definition
 of a blender as a hyperbolic basic set of a diffeomorphism
 which has the `thick' stable (or unstable) set
 (see also \cite[Definition 6.11]{BDV}).
In \cite{BBD16}, Bochi, Bonatti, and D\'iaz
 also provided more abstract definitions of blenders
 from geometrical and dynamical view points.
Let us recall that the definition of geometric blenders in \cite{BBD16}.
A hyperbolic set $\Lambda$ of a diffeomorphism $f$
 is called a {\it geometric $cu$-blender} of the $uu$-index $d^{uu}$ if
 its unstable index is greater than $d^{uu}$
 and there exists a $C^1$-open family $\cD$ of $d^{uu}$-dimensional disks
 and a $C^1$-neighborhood $\cU$ of $f$
 such that each disks in $\cD$ intersects with
 the stable set of the continuation of $\Lambda$ for any $g \in \cU$.
We call the open family $\cD$ of disks the {\it superposition region}.
Since the codimension of the stable manifold of each point of $\Lambda$
 equals to the unstable index of $\Lambda$ and it is greater than $d^{uu}$,
 no intersection of each stable manifold of $\Lambda$
 and a $d^{uu}$-dimensional disk is transverse.
However,
 the stable set of a blender of $uu$-index $d^{uu}$
 intersects with a $d^{uu}$-dimensional disk in the superposition region
 persistently, just like a submanifold of codimension $d^{uu}$
 transverse to the $d^{uu}$-dimensional disk.
This property enables invariant manifolds of lower dimensions to
 have $C^1$-robust non-transverse intersection,
 which is a source of $C^1$-robust phenomena mentioned above.

The aim of this paper is to give examples
 which exhibit the following robust phenomena
 with the help of mechanism inducing blenders:
\begin{enumerate}
 \item A pair of two regular Cantor sets in higher dimensions
 which exhibits $C^1$-stable intersection.
 \item A hyperbolic basic set which exhibits
 $C^2$-robust homoclinic tangency of the largest codimension.
\end{enumerate}
The former shows that Moreira's non-existence theorem
 of $C^1$-stable intersection of regular Cantor sets
 in the real line \cite{Mo11} does not hold in higher dimension.
The latter solves a problem on homoclinic tangency
 of large codimension posed by Barrientos and Raibekas \cite{BR18}.
In the next two subsections, we explain the background behind the examples
 and give the precise statement of the results.
In the last subsection of this section,
 we give an outline of the proof of the main theorems.

\subsection{Regular Cantor sets and their stable intersection}
\label{ss:Cantor}
A {\it Cantor set} is a topological space
 which is compact, totally disconnected, and without isolated points.
A famous example is the $1/3$-middle Cantor set
\begin{equation*}
 \Lambda_{1/3}=\left\{0.s_1s_2s_3\cdots=\sum_{n=1}^\infty 3^{-n}s_n \mid
 s_n \in \{0,2\} \right\},
\end{equation*}
 the set of real numbers in $[0,1]$ such that
 any digit of the ternary expansion is $0$ or $2$.
The Cantor set $\Lambda_{1/3}$ can be described as an invariant set of
 a dynamical system.
Take a $C^\infty$ map $f:\RR \ra \RR$
 such that $f(x)=3x$ if $x \in [0,1/3]$ and $f(x)=3x-2$ if $x \in [2/3,1]$.
Then, one can see that
\begin{equation*}
\Lambda_{1/3}=\bigcap_{n \geq 0}f^{-n}
\left(\left[0,\frac{1}{3}\right] \cup \left[\frac{2}{3},1\right]\right).
\end{equation*}
This dynamical construction of a Cantor set can be generalized as follows:
\begin{expl}
\label{expl:regular Cantor}
Let $(K_i)_{i \in I}$ be a family of
 mutually disjoint non-empty compact subsets
 in $\RR^m$ indexed by a finite set $I$,
 and a $C^1$ map $f$ from an open neighborhood 
 of the union $K=\bigcup_{i \in I}K_i$ to $\RR^m$
 such that $K \subset f(K_i) $ for any $i \in I$,
 and the restriction of $f$ to each $K_i$ is uniformly expanding,
 {\it i.e.}, there exists $\lambda_i>1$ such that
 $d(f(x),f(y)) \geq \lambda_i d(x,y)$ for any $x,y \in K_i$,
 where $d$ is a complete metric on $\RR^m$.
Then, the maximal forward invariant set
 $\Lambda^s(K,f)=\bigcap_{n \geq 0}f^{-n}(K)$ of $f$ in $K$
 is a Cantor set\footnote{From the viewpoint of fractal geometry,
 the family $((f|_{K_i})^{-1}:K \ra K_i)_{i \in I}$ of
 the inverse branches of $f$ is a {\it $C^r$ iterated function system}
 satisfying {\it the strong separation condition}.}.
To see this fact, let $\Sigma_I$ be the space of
 $I$-valued infinite sequences $(i_n)_{n \geq 0}$,
 with the product topology as a infinite product of the discrete space $I$.
This space is a Cantor set.
By uniform expansion on each $K_i$,
 the intersection $\bigcap_{n \geq 0}f^{-n}(K_{i_n})$
 contains a unique point $h(s)$ for each $s=(i_n)_{n \in 0} \in \Sigma^I$.
One can see that the map $h:\Sigma_I \ra \Lambda^s(K,f)$ is a homeomorphism.
\end{expl}
The Cantor set $\Lambda^s(K,f)$ is
 an example of a $C^r$-regular Cantor set\footnote{
In some literature ({\it e.g.} \cite{PT93}),
 a regular Cantor set is called a {\it dynamically defined} Cantor set.}.
We call a Cantor set $\Lambda$ in $\RR^m$ {\it $C^r$-regular}
 if there exist a pair $((K_i)_{i \in I},f)$
 of family of compact subsets of $\RR^m$ 
 and a $C^r$ locally expanding map $f$ as above,
 and $a_{ij} \in \{0,1\}$ for each $i,j$ such that
\begin{equation*}
\Lambda=\{h(s) \mid s=(i_n)_{n \geq 0} \in \Sigma^I, a_{i_{n+1}i_n}=1
 \text{ for all } n \geq 0\},
\end{equation*} 
 where $h:\Sigma_I \ra \Lambda^s(K,f)$ is the homeomorphism above.
Remark that $\Lambda=\Lambda^s(K,f)$ if all $a_{ij}=1$.

We define the {\it arithmetic sum} $K_1+K_2$
 and the {\it arithmetic difference} $K_1-K_2$
 of Cantor sets $K_1,K_2$ in the real line $\RR$ by
\begin{align*}
 K_1 \pm K_2=\{x \pm y \mid x \in K_1, y \in K_2\}.
\end{align*}
The arithmetic sum (or difference) of Cantor sets naturally
 appears in several mathematical problems,
 for instance, Markov and Lagrange spectra
 in the number theory  ({\it e.g.}, \cite{Mo98}),
 homoclinic bifurcation in the theory of smooth dynamical systems
 ({\it e.g.}, \cite{Mo96,Mo98}),
 and the spectrum of a Hamiltonian related to quasi-crystals
 ({\it e.g.}, \cite{DGS15}). 
The difference $K_1+K_2$ contains a real number $t$
 if and only if two Cantor sets $K_1$ and
 $K_2+t=\{y+t \in y \in K_2\}$ intersect.
Thus, a problem on the arithmetic difference (or sum) of Cantor sets
 can be easily translated into a problem on intersection of Cantor sets.

We can define $C^r$-perturbation of a regular Cantor set
 by $C^r$-perturbation of the locally expanding map
 that defines the Cantor set.
We say that a pair of regular Cantor sets has
 {\it $C^r$-stable intersection}
 if any $C^r$-perturbation of the pair has non-empty intersection
 (see Section \ref{sec:Cantor} for the precise definition).
If regular Cantor sets $K_1$ and $K_2$ in $\RR$
 have $C^r$-stable intersection for some $r\geq 1$,
 then $K_1$ and $K_2+t$ intersect for any $t$ sufficiently close to zero,
 and hence, the arithmetic difference $K_1-K_2$ contains
 a neighborhood of the origin.

In bifurcation theory, stability of intersection of
 regular Cantor sets plays an important role.
In \cite{Ne74} (see also \cite{Ne79,PT93}),
 Newhouse defined a numerical invariant called {\it thickness}
 for Cantor sets in $\RR$
 and proved that a pair of $C^2$-regular Cantor sets with large thickness
 has $C^2$-stable intersection.
He applied this result to show persistence of homoclinic tangency
 and abundance of diffeomorphisms with
 infinitely many attracting periodic orbits
 for $C^2$ surface diffeomorphisms.
His result was generalized to higher dimensional case
 by Palis and Viana \cite{PV94}
 and by Gonchenko, Turaev, and Shilnikov
 \cite{GTS93-1,GTS93-2}.
Kiriki and Soma \cite{KS12} applied
 Newhouse's thickness criterion to prove
 persistence of heterodimensional tangency.
In \cite{Bu93,Bu97}, Buzzard studied stable intersection
 of regular Cantor sets generated by holomorphic maps
 of the complex line $\mathbb{C}$
 and proved results analogous to Newhouse's one in holomorphic setting
\footnote{We refer \cite{AM-pre,Bi-pre} for recent progress
 in holomorphic case.}.

Ures \cite{Ur95} showed that
 thickness of $C^1$ generic regular Cantor sets in the real line
 is zero.
This means that Newhouse's thickness criterion 
 is useless for finding a pair having $C^1$-stable intersection.
Moreira proved the following negative result on $C^1$-stable intersection
 of Cantor sets in the real line.
\begin{thm}
[{Moreira \cite{Mo11}}]
\label{thm:Moreira}
There exists no pair of regular Cantor sets in $\RR$
 which has $C^1$-stable intersection.
\end{thm}
The first main result of this paper asserts that
 an analog of Theorem \ref{thm:Moreira}
 {\it does not} hold in higher dimension.
\begin{thmA}
\label{Mthm:Cantor}
For any positive integers $l,m$ and any positive real number $\delta$,
 there exists a pair $(\Lambda_1,\Lambda_2)$
 of regular Cantor sets in $\RR^{l+m}$
 which has $C^1$-stable intersection
 and satisfies $\udim_B(\Lambda_1)<l+\delta$, $\udim_B(\Lambda_2)<m+\delta$,
 where $\udim_B(K)$ is {\it the upper box dimension}
\footnote{
We refer \cite[Section 2.1]{Ba}
 for the definition and basic properties of
 the upper box dimension.
 It is known
 that the upper box dimension coincides with the Hausdorff dimension
 for any $C^2$ regular Cantor set in $\RR$
 (and more general `conformal' case).
 See \cite[Theorem 4.1.7]{Ba} for the proof.}
 of a compact subset $K$ of $\RR^m$.
\end{thmA}
We remark on a lower bound of the upper box dimensions
 for stable intersection.
For compact subsets $K_1$ and $K_2$
 in the $(l+m)$-dimensional Euclidean space
 $\RR^{l+m}$, it is easy to see that the difference $K_1-K_2$
 satisfies $\udim_B(K_1-K_2)\leq \udim_B(K_1)+\udim_B(K_2)$.
This means that if  $\udim_B(K_1)+\udim_B(K_2)<l+m$ then $K_1-K_2$
 has empty interior,
 and hence, $K_1$ and the translation $K_2+t=\{y+t \mid y \in K_2\}$
 of $K_2$ do not intersect for dense $t \in \RR^{l+m}$.
Therefore, any pair $(K_1,K_2)$ having stable intersection
 must satisfy the inequality $\udim_B(K_1)+\udim_B(K_2) \geq l+m$.
In particular,
 the upper box dimensions of Cantor sets in Theorem A
 are almost best possible.

For regular Cantor sets in the real line,
 Moreira and Yoccoz proved abundance of pairs
 having $C^2$-stable intersection.
\begin{thm}
[{Moreira and Yoccoz \cite{MY01}}]
Let $\mc{K}_1$ be the space of pairs
 $(\Lambda_1,\Lambda_2)$ of $C^2$-regular Cantor sets
 in the real line such that $\udim_B(\Lambda_1)+\udim_B(\Lambda_2)>1$.
Then, there exists an open and dense subset $\cU$ of $\mc{K}_1$ such that
\begin{equation*}
 I_S(\Lambda_1,\Lambda_2)=\{t \in \RR \mid
 \text{$\Lambda_1$ and $\Lambda_2+t$ have $C^2$-stable intersection}\}
\end{equation*}
 is a dense subset of  $\Lambda_1-\Lambda_2$
 for any $(\Lambda_1,\Lambda_2)$ in $\cU$.
\end{thm}
It is natural to ask whether the analogy for higher dimension holds or not.
\begin{qtn}
For $r \geq 1$,
 let $\mc{K}_m^r$ be the space of pairs $(\Lambda_1,\Lambda_2)$
 of $C^r$-regular Cantor sets in $\RR^m$
 such that $\udim_B(\Lambda_1)+\udim_B(\Lambda_2)>m$.
Does there exist an open and dense subset $\cU$ of $\mc{K}_m^r$
 such that $I_S(\Lambda_1,\Lambda_2)$ is a dense subset
 of $\Lambda_1-\Lambda_2$
 for any pair $(\Lambda_1,\Lambda_2)$ in $\cU$?
\end{qtn}
The case $r=1$ might be the most interesting case.

\subsection{Robust homoclinic tangency of large codimension}
\label{ss:tangency}
Let $f$ be a $C^r$ diffeomorphism of a manifold $M$
 and $\Lambda$ a topologically transitive hyperbolic set
 of unstable index $d^u$
 (see Sections \ref{sec:maximal} for the precise definition).
We say that $\Lambda$ exhibits {\it homoclinic tangency} at $q \in M$
 if there exist $p_1,p_2 \in \Lambda$ such that
 the stable manifold $W^u(p_1)$ and the unstable manifold $W^s(p_2)$
 intersect at $q$ non-transversely,
 {\it i.e.}, $T_q M \neq T_q W^u(p_1)+ T_q W^s(p_2)$.
When $\Lambda$ exhibits homoclinic tangency at $q$,
 we define {\it the codimension $c(q)$ of tangency} at $q$ by
\begin{align*}
 c(q) & = \dim M -\dim (T_q W^u(p_1)+ T_q W^s(p_2)).
\end{align*}
Since $\dim T_q W^u(p_1)=d^u$ and $\dim T_q W^s(p_2)=\dim M-d^u$,
 the inequality
\begin{equation*}
 c(q) \leq \min\{d^u,\dim M-d^u\} \leq \frac{1}{2}\dim M
\end{equation*}
 holds.
The codimension of tangency quantifies
 how far from transverse the intersection of
 the stable and unstable manifolds is.
Remark that $c(q)=\dim M/2$ if and only if $T_q W^u(p_1)=T_q W^s(p_2)$.
We say that $\Lambda$ exhibits homoclinic tangency of codimension $d$
 if there exists a point $q \in M$
 where $\Lambda$ exhibits homoclinic tangency with $c(q)=d$.

Suppose that a hyperbolic invariant set $\Lambda$ of $f$
 is locally maximal {\it i.e}, $\Lambda=\bigcap_{n \in \ZZ}f^{-n}(U)$
 for some neighborhood $U$ of $\Lambda$.
It is known that there exists a $C^r$-neighborhood $\cU$ of $f$
 such that the {\it continuation} $\Lambda(g)=\bigcap_{n \in \ZZ}g^{-n}(U)$
 of $\Lambda$ at $g$ is a hyperbolic set for any $g \in \cU$.
We say that the hyperbolic set $\Lambda$ exhibits
 {\it $C^r$-robust homoclinic tangency} of codimension $d$
 if the continuation $\Lambda(g)$
 exhibits homoclinic tangency of codimension $d$
 for any diffeomorphism $g$ sufficiently $C^r$-close to $f$.
Study of homoclinic tangency of higher codimension was initiated
 by Barrientos and Raibekas in \cite{BR17}.
They constructed a diffeomorphism
 which exhibits $C^2$-robust homoclinic tangency of codimension $c_T$
 for any manifold $M$ with $\dim M \geq 4$
 and any $1 \leq c_T \leq \lfloor \frac{1}{2}\dim M \rfloor-1$,
 where $\lfloor x \rfloor$ is the maximal integer not greater than $x$.
Barrientos and Raibekas
 also gave diffeomorphisms exhibiting $C^2$-robust homoclinic tangency
 of large codimension in symplectic setting in \cite{BR18}
 and outside partially hyperbolic region in \cite{BR20}.
Buzzi, Crovisier, and Fisher \cite{BCF18} 
 and Catalan \cite{Ca19} used
 (non-robust) homoclinic tangency of largest codimension
 to obtain a lower estimate of topological entropy
 for $C^1$ generic diffeomorphisms far from partially hyperbolic ones.

The following question is quite natural 
 since a manifold having a diffeomorphism
 with homoclinic tangency of codimension two
 must be at least four-dimensional.
\begin{qtn}
[Barrientos and Raibekas {\cite[p.4369]{BR17},
 see also \cite[Question 3]{BP-pre}}]
\label{qtn:BR}
Does there exist a diffeomorphism on a four-dimensional manifold
 which exhibits
 $C^2$-robust homoclinic tangency of codimension two?
\end{qtn}
More generally, one may ask
\begin{qtn}
Does there exist a diffeomorphism exhibiting
 $C^2$-robust homoclinic tangency of the largest codimension,
 {\it i.e.}, $\lfloor \frac{m}{2} \rfloor$ for
 a $m$-dimensional manifold?
\end{qtn}
The second main theorem of this paper solves this question.
\begin{thmB}
For any manifold $M$ with $\dim M \geq 4$,
 there exists a $C^\infty$ diffeomorphism which exhibits
 $C^2$-robust homoclinic tangency of codimension
 $\lfloor \frac{1}{2} \dim M \rfloor$.
\end{thmB}

\subsection{Outlines of Proofs}
Suppose that two diffeomorphisms $f_1$ and $f_2$ of $\RR^{l+m}$
 admit partially hyperbolic $cu$-blenders $\Lambda_1$ and $\Lambda_2$
 of $uu$-indices $m$ and $l$ respectively such that
\begin{itemize}
 \item the center-stable direction of $\Lambda_1$ is horizontal,
 the strong unstable direction of $\Lambda_1$ is vertical,
 \item the center-stable direction of $\Lambda_2$ is vertical, and
 the strong unstable direction of $\Lambda_2$ is horizontal
\end{itemize}
 with respect to the splitting $\RR^{l+m}=\RR^l \times \RR^m$.
An example of such a pair $(f_1,f_2)$ for the case $l=m=2$
 is given by the product $f_1(x_1,x_2,x_3,x_4)=(f(x_1,x_2,x_3),9x_4)$
 of a diffeomorphism $f$ in Example \ref{expl:blender horseshoe}
 and a linear expanding map $g(x_4)=9x_4$,
 and the transpose $f_2=T \circ f_1 \circ T^{-1}$ of $f_1$
 by the involution $T(x_1,x_2,x_3,x_4)=(x_3,x_4,x_1,x_2)$
 associate with the splitting $\RR^4=\RR^2 \times \RR^2$.
As we mentioned above,
 the stable set of a $cu$-blender behaves like a manifold
 of lower codimension virtually. 
In this case, the stable sets $W^s(\Lambda_1,f_1)$ and
 $W^s(\Lambda_2,f_2)$
 behave like a horizontal $l$-dimensional manifold
 and a vertical $m$-dimensional manifold respectively.
Although the stable sets are not a manifold and
 their topological dimensions are smaller than $l$ and $m$,
 one might expect that an analog of
 $C^1$-persistence of transverse intersection of
 horizontal and vertical manifolds
 holds for the stable sets $W^s(\Lambda_1,f_1)$ and
 $W^s(\Lambda_2,f_2)$.
We realize this naive idea under a suitable setting
 in Theorem \ref{thm:intersection}, which is a keystone
 to prove Theorems A and B.
In order to describe the setting,
 we introduce a {\it blending machine},
 which is a system of rectangles with respect to the splitting
 which induces a $cu$-blender whose stable set is horizontal\footnote{
 Describing a blender by rectangles dates back to
 the original work by Bonatti and D\'iaz in \cite{BD96}.
In \cite{BR17}, Barrientos and Raibekas also used a system of rectangles
 to describe a blender.}. 
Theorem \ref{thm:intersection} asserts 
 $C^1$-stability of the intersection of the forward invariant sets of
 a horizontal blending machine and a vertical bending machine
 like transverse intersection of horizontal and vertical manifolds.
We define blending machines
 and prove Theorem \ref{thm:intersection} in Section \ref{sec:intersection}.

We prove Theorem A
 by applying Theorem \ref{thm:intersection}
 to locally expanding maps which generate regular Cantor sets
 as their forward invariant sets.
This will be done in Section \ref{sec:Cantor}.

To prove Theorem B,  we apply Theorem \ref{thm:intersection}
 to the lifts of certain hyperbolic diffeomorphisms
 to the Grassmannian bundle, following
 the strategy of Barrientos and Raibekas in \cite{BR17},
 where they gave a hyperbolic set which exhibits
 $C^2$-robust homoclinic tangency of higher codimension as mentioned above.
A simple but useful observation by Barrientos and Raibekas
 was that intersection of the lifts of the stable
 and unstable sets of a hyperbolic set to the Grassmannian bundle
 induces homoclinic tangency of higher codimension\footnote{
In \cite{Ar02}, Arai used the same observation for the projective bundle
 to detect homoclinic tangency of codimension one.}.
Barrientos and Raibekas constructed a hyperbolic set
 whose lift to the Grassmannian bundle contains a $cu$-blender.
This blender in the Grassmannian bundle
 provides $C^1$-stable intersection of the lift of the stable set
 and the lift of the unstable manifold of a point in the hyperbolic set.
By the observation above,
 this implies $C^2$-robust homoclinic tangency of higher codimension.
For the largest codimension case,
 it seems impossible to construct a hyperbolic set
 such that the lift contains a blender
 whose virtual dimension as a manifold-like set is sufficiently large
 for having robust intersection with the lift of the unstable manifold
 of a point.
This is the reason why we can not apply Barrientos and Raibekas' proof
 to the largest codimension case directly.
In this paper, we overcome this difficulty
 with the help of Theorem \ref{thm:intersection}.
Although the dimension of the lift of the unstable manifold
 is not large enough,
 the virtual dimension of the unstable set of a $cs$-blender
 in the Grassmannian bundle might be large enough as a manifold-like set
 so that the sum of the virtual dimensions of the stable and the unstable
 sets of two blenders equal to the dimension of the manifold.
In fact, if one can find a pair of $cu$- and $cs$-blenders
 in the Grassmaniann bundle
 such that the stable set of a $cu$-blender is horizontal
 and the unstable set of $cs$-blender is vertical,
 then they intersects $C^1$-stable by Theorem \ref{thm:intersection}.
Then, $C^2$-robust homoclinic tangency of the largest codimension
 follows from Barrientos and Raibekas' observation above.
We prove Theorem B in Section \ref{sec:robust} along this strategy.

\subsection{Acknowledgment}
This paper was partially supported by the JSPS Kakenki Grants 18K03276.
The author would like to thank P.G.Barrientos and A. Raibekas,
 who explained their construction of a blender in the Grassmannian bundle
 and let me know their question when they visited Kyoto in December 2018.
The author is also grateful to an anonymous referee for many suggestions
 and comments which improved the article greatly.

\section{Stable intersection of the forward invariant sets}
\label{sec:intersection}

\subsection{Splittings and the cone condition}
\label{sec:cones}

For manifolds $M_1$ and $M_2$ and $r \geq 1$,
 let $C^r(M_1,M_2)$ be the set of $C^r$ maps from $M_1$ to $M_2$
 with the compact-open $C^r$-topology.
For a subset $S$ of $M_1$, we say that a map $f:S \ra M_2$ is $C^r$
 if it extends to a $C^r$ map from an open neighborhood of $S$ to $M_2$.
For $m \geq 1$, we denote the $m$-dimensional Euclidean space by $\RR^m$
 and {\bf the box norm} on $\RR^m$ by $\|\cdot\|$,
 {\it i.e.}, $\|(x_1,\dots,x_m)\|=\max\{|x_1|,\dots,|x_m|\}$.
We identify the tangent space $T_x \RR^m$ at each $x \in \RR^m$
 with $\RR^m$ in a natural way.

Fix positive integers $l,m$.
We say that a pair $(P,Q)$ is an {\it $(l,m)$-splitting} of
 an $(l+m)$-dimensional open manifold $M$
 if $P$ and $Q$ are $C^1$ maps from $M$ to $\RR^l$ and $\RR^m$ respectively,
 and the map $(P,Q): M \ra \RR^{l+m}$
 given by $(P,Q)(x)=(P(x),Q(x))$
 is a diffeomorphism onto an open subset of $\RR^{l+m}$.
Let $M$ be a manifold with an $(l,m)$-splitting $(P,Q)$.
For $x \in M$ and $\theta>0$,
 we define the {\it $\theta$-cone} $C(x,\theta,P,Q)$ by
\begin{align*}
 C(x,\theta,P,Q)
  & = \{v \in T_x M \mid \|DQ_x v\| \leq \theta \|DP_x v\|\}.
\end{align*}
Let $M_1$, $M_2$ be manifolds with $(l,m)$-splitting $(P_1,Q_1)$,
 $(P_2,Q_2)$ respectively,
 $U$ an open subset of $M_1$, and $f:U \ra M_2$ a $C^1$ embedding.
For $\theta>0$ and a subset $S$ of $U$,
 we say that $f$ satisfies the $\theta$-cone condition on $S$ 
 with respect to the splittings $(P_1,Q_1)$ and $(P_2,Q_2)$
 if there exists $0<\theta'<\theta$ such that
\begin{gather*}
 Df^{-1}(C(f(x),\theta,P_2,Q_2)) \subset C(x,\theta',P_1,Q_1), \\
 Df(C(x,\theta,Q_1,P_1)) \subset C(f(x),\theta',Q_2,P_2)
\end{gather*}
 for any $x \in S$.
When $M_1=M_2$ and $(P_1,Q_1)=(P_2,Q_2)$,
 we just say that $f$ satisfies the $\theta$-cone condition on $S$ 
 with respect to $(P_1,Q_1)$.
For $\theta,\lambda,\mu>0$,
 we say that $f$ satisfies the $(\theta,\lambda,\mu)$-cone condition
 if it satisfies the $\theta$-cone condition on $S$
 and there exists $\epsilon>0$ such that
\begin{alignat*}{4}
\|D(P_1 \circ f^{-1}) v \| & \geq (\lambda+\epsilon) \|DP_1 v\|, & \hsp
\|D(Q_2 \circ f) w \| & \geq (\mu+\epsilon) \|DQ_1 w \|
\end{alignat*}
 for any $x \in S$,
 $v \in C(f(x),\theta,P_2,Q_2)$, and $w \in C(x,\theta,Q_1,P_1)$.
\begin{rmk}
If a $C^1$ embedding $f:U \ra M_2$ satisfies
 the $(\theta,\lambda,\mu)$-cone condition on $S$
 with respect to $(P_1,Q_1)$ and $(P_2,Q_2)$,
 then the inverse $f^{-1}:f(U) \ra M_1$ satisfies 
 the $(\theta,\mu,\lambda)$-cone condition on $f(S)$
 with respect to $(Q_2,P_2)$ and $(Q_1,P_1)$.
\end{rmk}
\begin{rmk}
If $f$ satisfies the $(\theta,\lambda,\mu)$-cone condition
 on a compact subset $K$ of $U$,
 then there exist a neighborhood $U_K$ of $K$
 and $C^1$-neighborhood $\cU$ of $f$ such that
 any $g \in \cU$ satisfies
 the $(\theta,\lambda,\mu)$-cone condition on $U_K$.
\end{rmk}
\begin{rmk}
\label{rmk:cone hyp}
Let $M$ be a manifold with an $(l,m)$-splitting $(P,Q)$
 and $U$ an open subset of $M$.
If a $C^1$-embedding $f:U \ra M$ satisfies the $(\theta,1,1)$-cone condition
 with respect to $(P,Q)$ on a compact $f$-invariant subset $\Lambda$ of $U$
 for some $\theta>0$, then there exist
 a $Df$-invariant continuous splitting $E^s \oplus E^u$ of $TM$ on $\Lambda$ 
 and constants $\kappa>1$, $0<\gamma<1$,
 which satisfy the following conditions for any $x \in \Lambda$:
\begin{itemize}
 \item $\dim E^s(x)=l$, $\dim E^u(x)=m$,
 \item $E^s(x) \subset C(x,\theta,Q,P)$, $E^u(x) \subset C(x,\theta,P,Q)$,
 \item $\|Df^n v\| \leq \kappa\gamma^n \|v\|$
 and $\|Df^{-n} w\| \leq \kappa\gamma^n \|w\|$
 for $v \in E^s(x)$, $w \in E^u(x)$ and $n \geq 0$.
\end{itemize}
In particular, $\Lambda$ is a {\it hyperbolic set}
 with the {\it the hyperbolic splitting} $E^s \oplus E^u$.
See Section \ref{sec:maximal} for the precise definition
 of a hyperbolic set and its hyperbolic splitting.
\end{rmk}

\subsection{Rectangles}
\label{sec:rectangles}
Fix positive integers $l,m$
 and an $(l+m)$-dimensional manifold $M$.
For $C^1$-maps $G:M \ra \RR^l$, $H:M \ra \RR^m$
 (the pair $(G,H)$ does not need to be an $(l,m)$-splitting of $M$),
 we say that a subset $R$ of $M$ is a $(G,H)$-rectangle
 if $(G,H)(R)=G(R) \times H(R)$
 and the map $(G,H):R \ra G(R) \times H(R)$ is a diffeomorphism.
\begin{lemma}
\label{lemma:rectangle stability}
Let $R \subset M$ be a compact $(G,H)$-rectangle.
For any open neighborhood $U_R$ of $R$
 and any compact subset $K_R$ of $\Int R$,
 there exists a $C^1$-neighborhood $\cU$ of $(G,H)$
 such that any $(G',H') \in \cU$ admits
 a compact $(G',H')$-rectangle $R'$ satisfying 
 $G'(R')=G(R)$, $H'(R')=H(R)$, $K_R \subset \Int R'$, and $R' \subset U_R$.
\end{lemma}
\begin{proof}
If a pair $(G',H')$ is sufficiently $C^1$-close to $(G,H)$,
 then $(G',H')$ maps a neighborhood of $R$ to
 a neighborhood of $G(R) \times H(R)$ diffeomorphically.
Then, $R'=(G',H')^{-1}(G(R) \times H(R))$ is a $(G',H')$-rectangle
 such that $G'(R')=G(R)$ and $H'(R')=H(R)$.
If $(G',H')$ is sufficiently close to $(G,H)$,
 then $K_R \subset \Int R'$ and $R' \subset U_R$.
\end{proof}

Rectangles 'skewed' by a diffeomorphism and their crossing
 are important in this paper.
\begin{figure}
\begin{center}
 \includegraphics[scale=0.6]{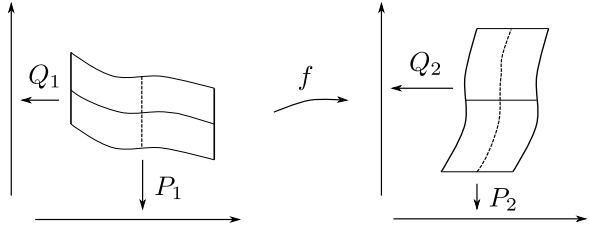} 
\end{center}
\caption{A $(P_1,Q_2 \circ f)$-rectangle -- a rectangle skewed by $f$}
\end{figure}
Fix a manifold $M_n$ with an $(l,m)$-splitting $(P_n,Q_n)$
 for $n=1,\dots,N+1$
 and a $C^1$ embedding $f_n:U_n \ra M_{n+1}$ for $n=1,\dots,N$,
 where $U_n$ is an open subset of $M_n$.
Define a map $F_n:U_{1,n} \ra M_{n+1}$
 by $F_n=f_n \circ \dots \circ f_1$ for $n=1, \dots, N$,
 where $U_{1,n}$ is the maximal open subset of $U_1$
 on which the composition $f_n \circ \dots \circ f_1$ is well-defined.
Put $U_{1,0}=U_1$ and let $F_0:U_1 \ra M_1$ be the inclusion map.
We give two lemmas for such sequence of maps.
\begin{lemma}
\label{lemma:rectangle composition}
Let $(R_n)_{n=1}^N$ be a sequence
 such that $R_n$ is a compact $(P_n,Q_{n+1} \circ f_n)$-rectangle in $U_n$
 and
\begin{alignat*}{4}
P_{n+1}(f_n(R_n)) & \subset \Int P_{n+1}(R_{n+1}), & \quad
Q_{n+1}(R_{n+1}) & \subset \Int Q_{n+1}(f_n(R_n)),
\end{alignat*}
 for $i=1,\dots, N$.
Suppose that there exists a constant $0<\theta \leq 1$
 such that $f_n$ satisfies the $\theta$-cone condition
 with respect to $(P_n,Q_n)$ and $(P_{n+1},Q_{n+1})$ on $R_n$
 for any $n=1,\dots,N$.
Then, $R_*=\bigcap_{n=1}^N (F_{n-1})^{-1}(R_n)$ is
 a $(P_1, Q_{N+1} \circ F_N)$-rectangle such that
\begin{alignat*}{4}
 P_1(R_*)
 & = P_1(R_1), & \quad
 Q_{N+1}(F_N(R_*)) & = Q_{N+1}(f_N(R_N)).
\end{alignat*}
\end{lemma}
\begin{proof}
Proof is done by induction of $N$.
The case $N=1$ is trivial.
Suppose that the lemma holds for the case $N-1$.
Then, $R'_*=\bigcap_{n=1}^{N-1}(F_{n-1})^{-1}(R_n)$
 is a $(P_1,Q_N \circ F_{N-1})$-rectangle
 such that $P_1(R'_*)=P_1(R_1)$ and
 $Q_N(F_{N-1}(R'_*))=Q_N(f_{N-1}(R_{N-1}))$.
Since each $f_n$ satisfies the $\theta$-cone condition
 with respect to $(P_n,Q_n)$ and $(P_{n+1},Q_{n+1})$
 on the compact subset $R_n$,
 there exists $0<\theta'<\theta$ such that
 $f_n$ satisfies the $\theta'$-cone condition on $R_n$ 
 with respect to the same splittings for each $n=1,\dots,N$.
Put $R_*=\bigcap_{n=1}^N (F_{n-1})^{-1}(R_n)$.
Then, $R_*=R_*' \cap (F_{N-1})^{-1}(R_N)$
 and $F_N$ satisfies the $\theta'$-cone condition on $R_*$
 with respect to the splittings $(P_1,Q_1)$ and $(P_{N+1},Q_{N+1})$.

First,
 we show that the restriction of $(P_1,Q_{N+1} \circ F_N)$ to $R_*$
 is an immersion.
Take $x \in R_*$ and $v \in T_x M$ with $D(P_1,Q_{N+1} \circ F_N)(v)=0$
 and we show $v=0$.
The vector $v$ satisfies that $DP_1v=0$ and $D(Q_{N+1} \circ F_N)v=0$.
The former implies that $v$ is contained in $C(x,\theta',Q_1,P_1)$.
By the $\theta'$-cone condition on $R_*$ for $F_N$,
 $DF_N v$ is contained in $C(F_N(x),\theta',Q_N,P_N)$,
 that is,
 $\|D(P_{N+1} \circ F_N) v\| \leq \theta'\|D(Q_{N+1} \circ F_N) v\|$.
Since $D(Q_{N+1} \circ F_N) v=0$, we have
 $D(P_{N+1} \circ F_N) v=0$, and hence, $D((P_{N+1},Q_{N+1}) \circ F_N)v=0$.
This implies $v=0$ since $(P_{N+1},Q_{N+1})$ and $F_N$ are embeddings.
Therefore, 
 the kernel of $D(P_1,Q_{N+1} \circ F_N)_x$ is trivial for any $x \in R_*$,
 and hence, the map
 $(P_1,Q_{N+1} \circ F_N):R_* \ra P_1(R_1) \times (Q_{N+1} \circ f_N)(R_N)$
 is an immersion.
\begin{figure}
\begin{center}
 \includegraphics[scale=0.6]{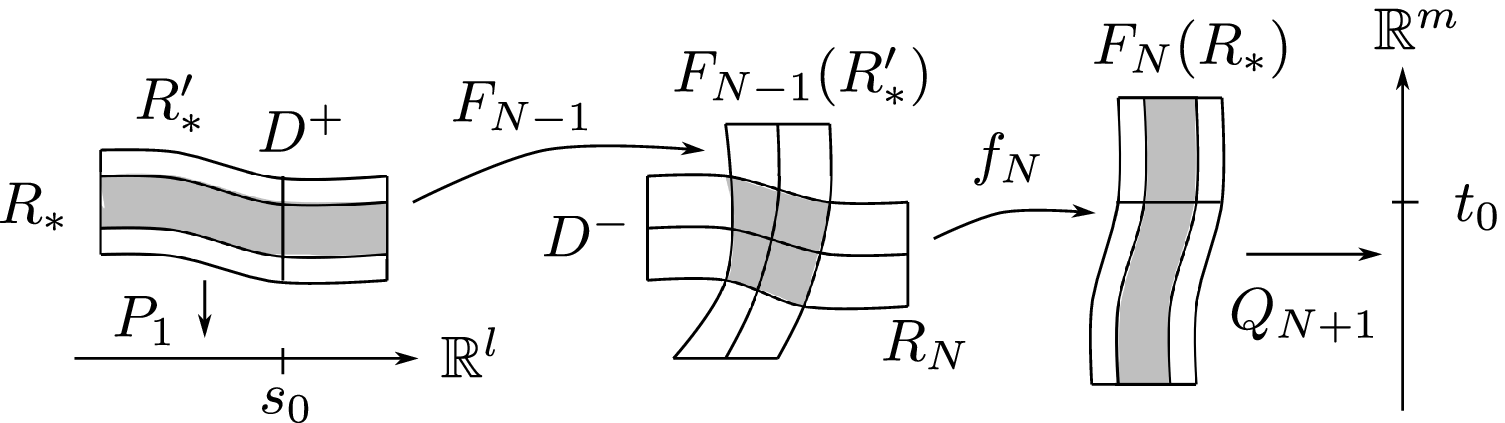} 
\end{center}
 \caption{Proof of Lemma \ref{lemma:rectangle composition}}
\end{figure}

Fix $s_0 \in P_1(R_1)$ and $t_0 \in (Q_{N+1} \circ f_N)(R_N)$ and we put
\begin{align*}
D^+  & =\{x \in R_1 \mid P_1(x)=s_0\}, &
D^-  & =\{y \in R_N \mid (Q_{N+1} \circ f_N)(y)=t_0\}.
\end{align*}
We claim that $F_{N-1}(D^+)$ intersects with $D^-$ at a unique point.
The claim implies that
 $(P_1,Q_{N+1} \circ F_N)$ is a bijection from $R_*$
 to $P_1(R_1) \times (Q_{N+1} \circ f_N)(R_N)$.
Since the map is an immersion on a neighborhood of $R_*$,
 it is a diffeomorphism.
Therefore, the claim implies the lemma for the case $N$
 and the induction completes the proof of the lemma.

Let us show the claim.
Since $R_*'$ is a $(P_1, Q_N \circ F_{N-1})$-rectangle,
 $F_{N-1}(R_*')$ is a $(P_1 \circ F_{N-1}^{-1},Q_N)$-rectangle.
This implies that there exists a $C^1$ map
 $g^+_{s_0}:Q_N(F_{N-1}(R_*')) \ra P_N(F_{N-1}(R_*'))$ 
 such that $(P_N,Q_N)(F_{N-1}(D^+))$
 is the graph of $g^+_{s_0}$ on $Q_N(F_{N-1}(R_*'))$,
 {\it i.e.},
\begin{equation}
\label{eqn:composition 1}
(P_N,Q_N)(F_{N-1}(D^+))=\{(g^+_{s_0}(t),t) \mid t \in Q_N(F_{N-1}(R_*'))\}. 
\end{equation}
By the $\theta'$-cone condition for $F_{N-1}=f_{N-1} \circ \dots f_1$,
 we have $\|g^+_{s_0}\|\leq \theta'<\theta \leq 1$.
Similarly, there exists a $C^1$ function
 $g^-_{t_0}:P_N(R_N) \ra Q_N(R_N)$ such that
\begin{equation}
\label{eqn:composition 2}
(P_N,Q_N)(D^-)=\{(s,g^-_{t_0}(s)) \mid s \in P_N(R_N)\} 
\end{equation}
 and $\|g^-_{t_0}\|\leq \theta'<1$.
Since $P_N \circ F_{N-1}(R_*') \subset P_N(R_N)$
 and $Q_N(R_N) \subset Q_N \circ F_{N-1}(R_*')$,
 we can define a self-map $G$ on
 $P_N(R_N) \times (Q_N \circ F_{N-1})(R_*')$
 by $G(s,t)=(g^-_{s_0}(t),g^+_{t_0}(s))$.
This map is a uniform contraction,
 and hence, admits a unique fixed point $(s_*,t_*)$
 by the contracting mapping principle.
By (\ref{eqn:composition 1}) and (\ref{eqn:composition 2}),
 a point $x \in F_{N-1}(R_*') \cap R_N$
 is contained in $F_{N-1}(D^+) \cap D^-$
 if and only if $(P_N(x),Q_N(x))$ is a fixed point of $G$.
Therefore, $F_{N-1}(D^+)$ intersects with $D^-$ at the unique point
 $(P_N,Q_N)^{-1}(s_*,t_*)$.
This finishes the proof of the claim, and hence, of the lemma.
\end{proof}

For $m \geq 1$ and a subset $S$ of $\RR^m$,
 let $\diam S$ be the diameter of $S$ with respect to the metric
 induced by the box norm $\|\cdot\|$.
The following lemma gives a bound of the diameter of $f^n(R_*)$
 in the previous lemma.
\begin{lemma}
\label{lemma:rectangle size}
Let $R$ be a compact $(P_1,Q_{N+1} \circ F_N)$-rectangle in $U_{1,N}$.
Suppose that there exist constants $\theta,\lambda,\mu>0$
 such that $f_n$ satisfies
 the $(\theta,\lambda,\mu)$-cone condition on $F_{n-1}(R)$
 with respect to $(P_n,Q_n)$ and $(P_{n+1},Q_{n+1})$
 for each $n=1,\dots,N$.
Then,
\begin{align}
\label{eqn:size 1}
\diam Q_{n+1}(F_n(R))
 & \leq \theta \lambda^{-n} \diam P_1(R)
 + \mu^{-(N-n)}\diam Q_{N+1}(F_N(R)),\\
\label{eqn:size 2}
\diam P_{n+1}(F_n(R)) 
 & \leq  \lambda^{-n}\diam P_1(R) + \theta \mu^{-(N-n)}\diam Q_{N+1}(F_N(R))
\end{align}
 for any $0 \leq n \leq N$.
\end{lemma}
\begin{proof}
Take $x,x' \in R$.
Since $R$ is a $(P_1,Q_{N+1} \circ F_N)$-rectangle,
 there exists $x_* \in R$ such that
 $P_1(x_*)=P_1(x)$ and $Q_{N+1}\circ F_N(x_*)=Q_{N+1} \circ F_N(x')$.
Put
\begin{align*}
 R^- & = \{y \in R \mid Q_{N+1}(F_N(y))= Q_{N+1}(F_N(x_*))\}, \\
 R^+ & = \{y \in R \mid P_1(y)=P_1(x_*)\}.
\end{align*}
Then, $\{x,x_*\} \subset R^-$ and $\{x',x_*\} \subset R^+$.
By the $(\theta,\lambda,\mu)$-cone condition for $f_1,\dots,f_N$,
 the composition $f_{n'} \circ \dots \circ f_n$ satisfies the
 $(\theta,\lambda^{(n'-n)+1},\nu^{(n'-n)+1})$-cone condition
 for $1\leq n \leq n' \leq N$.
For any $y \in R^-$,
 any vector $v \in T_y R^-$ satisfies that $D(Q_{N+1} \circ F_N)(v)=0$,
 and hence, $DF_N v$ is contained in $C(F_N(y),\theta,P_{N+1},Q_{N+1})$.
By the cone condition, we have
\begin{align*}
 \|D(Q_{n+1} \circ F_n)v\|  & \leq \theta \|D(P_{n+1} \circ F_n) v\|
 \leq \theta\lambda^{-n} \|DP_1 v\|.
\end{align*}
 for any $0 \leq n \leq N$.
This implies that
\begin{align*}
 \diam Q_{n+1}(F_n(R^-)) & \leq \theta \lambda^{-n}\diam P_1(R^-). 
\end{align*}
For any $y' \in R^+$,
 any vector $v' \in T_{y'}R^+$ satisfies that $DP_1 v'=0$,
 and hence, it is contained in $C(y',\theta,Q_1,P_1)$.
By the cone condition, we have
\begin{align*}
 \|D(Q_{N+1} \circ F_N) w\|
  & \geq \mu^{N-n} \|D(Q_{n+1} \circ F_n) w\|
\end{align*}
 for any $0 \leq n \leq N$.
This implies that
\begin{align*}
 \diam Q_{n+1}(F_n(R^+)) & \leq \mu^{-(N-n)} \diam Q_{N+1}(F_N(R^+)).
\end{align*}
Since $\{x,x_*\} \subset R^-$ and $\{x',x_*\} \subset R^+$,
 we obtain that
\begin{align*}
 \|Q_{n+1} \circ F_n(x)-Q_{n+1} \circ F_n(x')\|
 & \leq \diam Q_{n+1}(F_n(R^-))+\diam Q_{n+1}(F_n(R^+)).
\end{align*}
Since $x$ and $x'$ are arbitrary points in $R$,
 we obtain (\ref{eqn:size 1}) by combining with the previous inequalities.
The inequality (\ref{eqn:size 2}) can be proved in the same way.
\end{proof}
The case $N=1$ of the above lemma implies
\begin{cor}
\label{cor:rectangle size}
Let $M_1, M_2$ be manifold with $(l,m)$-splittings $(P_1,Q_1), (P_2,Q_2)$
 respectively,
 $U$ an open subset of $M_1$, $f:U \ra M_2$ an $C^1$-embedding,
 and $R$ a $(P_1,Q_2 \circ f)$-rectangle in $U$.
Suppose that $f$ satisfies the $(\theta,\lambda,\mu)$-cone condition
 on $R$ with respect to $(P_1,Q_1)$ and $(P_2,Q_2)$.
Then,
\begin{align}
\label{eqn:size 1-1}
\diam Q_1(R) & \leq \theta \diam P_1(R)
 + \mu^{-1}\diam Q_2(f(R))\\
\label{eqn:size 2-1}
\diam P_2(f(R)) 
 & \leq \lambda^{-1}\diam P_1(R)+\theta\diam Q_2(f(R)).
\end{align}
\end{cor}

\subsection{The maximal Invariant sets in the union of rectangles}
\label{sec:maximal}
In this subsection,
 we give a simple criterion that 
 the maximal invariant set of a union of rectangles
 is a local maximal and topologically transitive hyperbolic set.

First, we recall some basic notions in topological dynamics.
Let $(X,d)$ be a metric space
 and $g:Y \ra X$
 a continuous map from an open subset $Y$ of $X$ to $X$.
We say that a subset $\Lambda$ of $Y$ is {\it $g$-invariant}
 if $g(\Lambda)=\Lambda$.
For a subset $S$ of $Y$, we define
 {\it the maximal invariant set} $\Lambda(S,g)$ and 
 {\it the maximal forward invariant set} $\Lambda^s(S,g)$ of $g$ in $S$ by
\begin{alignat*}{4}
 \Lambda(S,g) & =\bigcap_{n \in \ZZ}g^{-n}(S), & \hsp
 \Lambda^s(S,g) & =\bigcap_{n \geq 0}g^{-n}(S).
\end{alignat*}
Let $\Lambda$ be a $g$-invariant set.
 we define the stable set $W^s(\Lambda,g)$ of $\Lambda$ by
\begin{equation*}
 W^s(\Lambda,g)=\left\{ x \in U \mid
 d(g^n(x),\Lambda) \ra 0 \;(n \ra +\infty) \right\},
\end{equation*}
 where $d(x,\Lambda)=\inf_{y \in \Lambda}d(x,y)$ for $x \in X$.
We say that the invariant set $\Lambda$ is {\it locally maximal}
 if there exists an open neighborhood $U$ of $\Lambda$ such that
 $\Lambda=\bigcap_{n \in \ZZ}g^{-n}(U)$.
We say that $\Lambda$ is {\it topologically transitive} if
 any pair $(U,V)$ of non-empty open subsets of $\Lambda$
 admits an integer $n$ such that $g^n(U) \cap V \neq \emptyset$.
For any compact subset $K$ of $Y$,
 the stable set $W^s(\Lambda(K,g),g)$ contains $\Lambda^s(K,g)$
 since $\Lambda(K,g)=\bigcap_{n \geq 0}g^n(\Lambda^s(K,g))$.
The maximal $g$-invariant set $\Lambda(K,g)$ in $K$
 is locally maximal if and only if $\Lambda(K,g) \subset \Int K$.
When $g$ is a homeomorphism of $X$,
 we can define {\it the unstable set} $W^u(\Lambda,g)$
 of a $g$-invariant set $\Lambda$
 and {\it the maximal backward invariant set} $\Lambda^u(K,g)$
 of a compact set $K$ by
 $W^u(\Lambda,g)=W^s(\Lambda,g^{-1})$
 and $\Lambda^u(K,g)=\Lambda^s(K,g^{-1})$.

Next, we recall some basic notions in hyperbolic dynamics.
The standard references are \cite{KH} and \cite{Sh}.
Let $M$ be a smooth manifold and $f:U \ra M$
 a $C^r$ embedding from an open subset $U$ of $M$ to $M$.
We say that a compact invariant set $\Lambda$ of $f$ is {\it hyperbolic}
 if there exist a Riemannian metric on $M$, a constant $\lambda>1$,
 and a continuous splitting $TM|_\Lambda=E^s \oplus E^u$
 of the restriction of the tangent bundle $TM$ to $\Lambda$ such that
\begin{alignat*}{4}
 \|Df v\|_M &\leq \lambda^{-1}\|v\|_M, &\hsp
 \|Df w\|_M & \geq \lambda \|w\|_M
\end{alignat*}
 for any $n \geq 0$,
 $x \in \Lambda$, $v \in E^s(x)$, and $w \in E^u(x)$,
 where $\|\cdot \|_M$ is the norm associated with the Riemannian metric.
It is known that the dimensions of $E^u(x)$ and $E^s(x)$ are constant
 if $\Lambda$ is topologically transitive.
We call the dimensions of $E^s(x)$ and $E^u$
 {\it the stable index} and {\it the unstable index} of
 the hyperbolic set $\Lambda$.
Let $\Lambda$ be a compact hyperbolic invariant set of $f$.
For $x \in \Lambda$, we put
\begin{align*}
 W^s(x) & =\{y \in M \mid d(f^n(y),f^n(x)) \ra 0\; (n \ra +\infty)\},\\
 W^u(x) & =\{y \in M \mid d(f^n(y),f^n(x)) \ra 0\; (n \ra -\infty)\},
\end{align*}
 where $d$ is the distance induced from a Riemannian metric.
By the stable manifold theorem,
 $W^s(x)$ and $W^u(x)$ are injectively immersed submanifolds of $M$
 satisfying $T_x W^s(x)=E^s(x)$ and $T_x W^u(x)=E^u(x)$ for any $x$.
They are called
 {\it the stable manifold} and {\it the unstable  manifolds}
 of $f$ at $x$ respectively.
It is known that
 $W^s(\Lambda)  = \bigcup_{x \in \Lambda} W^s(x)$
 and $W^u(\Lambda)  = \bigcup_{x \in \Lambda} W^u(x)$
 if $\Lambda$ is a compact locally maximal hyperbolic set.

The following is a criterion that
 the maximal invariant set of the union of rectangles
 is a locally maximal and topologically transitive hyperbolic set.
\begin{prop}
\label{prop:transitive}
Let $M$ be a manifold with an $(l,m)$-splitting $(P,Q)$,
 $(K_i)_{i \in I}$ a family of mutually disjoint compact subsets of $M$
 indexed by a finite set $I$, and
 $f$ a $C^1$ embedding from an open neighborhood
 of the union $K=\bigcup_{i \in I} K_i$ to $M$.
Suppose that the map $f$ satisfies the $(\theta,1,1)$-cone condition on $K$
 for some $0<\theta \leq 1$,
 $K_i$ is $(P,Q \circ f)$-rectangles,
 $P(f(K_i))  \subset \Int P(K_j)$, and $Q(K_j) \subset \Int Q(f(K_i))$
 for any $i,j \in I$.
 then $\Lambda(K,f)$ is a topologically transitive
 hyperbolic set invariant set of unstable index $m$ for $f$
 with $\Lambda(K,f) \subset \Int K$.
\end{prop}
Remark that $\Lambda(K,f)$ is locally maximal
 since $\Int K$ is a neighborhood of $\Lambda(K,f)$ such that
 $\Lambda(K,f)=\bigcap_{n \in \ZZ}f^{-n}(\Int K)$.
\begin{proof}
[Proof of Proposition \ref{prop:transitive}]
As mentioned in Remark \ref{rmk:cone hyp},
 hyperbolicity of $\Lambda(K,f)$ follows from
 the $(\theta,1,1)$-cone condition on $K$.

We show that $\Lambda(K,f) \subset \Int K$.
By assumption, we have
\begin{align*}
 P(K_j \cap f(K_i))
 & = P(K_j) \cap P(f(K_i)) \subset \Int P(K_j),\\
 (Q \circ f)(K_j \cap f^{-1}(K_l))
 & = Q(f(K_j)) \cap Q(K_l) \subset \Int (Q \circ f)(K_j)
\end{align*}
 for any $i,j,l \in I$.
Since $K_j$ is a $(P,Q \circ f)$-rectangle,
 $(P, Q \circ f)$ is a diffeomorphism between $\Int K_j$ and
 $\Int P(K_j) \times \Int (Q \circ f)(K_j)$.
Hence, $f(K_i) \cap K_j \cap f^{-1}(K_l) \subset \Int K_j$.
This implies that $\Lambda(K,f) \subset \Int K$.
As we mentioned above,
 it follows that $\Lambda(K,f)$ is locally maximal.

We show that $\Lambda(K,f)$ is topologically transitive.
Let $I^{\ZZ}$ be the space of bi-infinite sequences $(i_n)_{n \in \ZZ}$
 of elements of $I$
 with the product topology of the discrete topological space $I$.
We define a map $h:\Lambda(K,f) \ra I^\ZZ$ by
 $h(x)_n=i$ if $f^n(x) \in K_i$.
It is easy to see that $h$ is continuous.
For a sequence $s=(i_n)_{i \in \ZZ}$ in $I^{\ZZ}$ and an integer $N \geq 0$,
 we put $K(s,N)=\bigcap_{|n| \leq N}f^{-n}(K_{i_n})$
 and $K(s)=\bigcap_{N \geq 0}K(s,N)=\bigcap_{n \in \ZZ}f^{-n}(K_{i_n})$.
Remark that $K(s)=h^{-1}(s)$.
Since $P(f(K_{i_n})) \subset \Int P(K_{i_{n+1}})$
 and $Q(K_{i_{n+1}}) \subset \Int Q(f(K_{i_n}))$ for any $n \in \ZZ$,
 Lemma \ref{lemma:rectangle composition} implies that
 $f^{-N}(K(s,N))$ is a $(P,Q \circ f^{2N})$-rectangle
 such that $P(f^{-N}(K(s,N)))=P(K_{i_{-N}})$
 and $Q(f^N(K(s,N)))=Q(K_{i_N})$ for any $N \geq 0$.
In particular,
 $(K(s,N))_{N \geq 0}$ is a nested sequence of non-empty compact sets,
 and hence, $K(s)=\bigcap_{N \geq 0}K(s,N)$ is non-empty.
This implies that the map $h$ is surjective.
Since $f$ satisfies the $(\theta,1,1)$-cone condition on $K$,
 it also satisfies the $(\theta,\nu,\nu)$-cone condition for some $\nu>1$.
Applying Lemma \ref{lemma:rectangle size}
 to the $(P,Q \circ f^{2N})$-rectangle $f^{-N}(K(s,N))$, we have
\begin{align*}
\diam P(K(s,N)) & \leq \nu^{-N}(\theta\diam P(K_{i_{-N}})+\diam Q(K_{i_N})),\\
\diam Q(K(s,N)) & \leq \nu^{-N}(\diam P(K_{i_{-N}})+\theta \diam Q(K_{i_N})),
\end{align*}
The right-hand sides of both inequalities are bounded
 by $\nu^{-N}(1+\theta)(\diam P(K)+\diam Q(K))$.
This implies that $\diam K(s,N)$ converges to zero as $N$ goes to infinity,
 and hence, $K(s)=h^{-1}(s)$ contains exactly one point.
Therefore, the map $h$ is bijective.
Since $h$ is a continuous bijective map between compact sets
 $\Lambda(K,f)$ and $I^\ZZ$, it is a homeomorphism.
It is easy to see that $h$ is a topological conjugacy
 between the restriction of $f$ to $\Lambda(K,f)$
 and the shift map on $I^\ZZ$.
It is well-known that the shift map on $I^\ZZ$ is topologically transitive.
Therefore, the restriction of $f$ to $\Lambda(K,f)$
 is topologically transitive.
\end{proof}

\subsection{Blending machine and a criterion to $C^1$-stable intersection}
\label{sec:blender}
Let $l,m$ be positive integers,
 $M$ an $(l+m)$-dimensional manifold with an $(l,m)$-splitting $(P,Q)$,
 $U$ an open subset of $M$,
 and $f:U \ra M$ a $C^1$ map.
We call a family $\mc{B}=(R_i)_{i \in I}$
 of compact subsets of $U$ indexed by a finite set $I$
 is a {\it blending machine} for $f$ with respect to $(P,Q)$ if
 there exist $\theta>0$ and  a $(P,Q)$-rectangle $Z$ such that
\begin{enumerate}
 \item the restriction of $f$ to an open neighborhood of $R_i$
 is a $C^1$ embedding which satisfies the $\theta$-cone condition
 with respect to $(P,Q)$ for any $i \in I$,
 \item  $R_i$ is a $(P,Q \circ f)$-rectangle such that
 $Q(R_i) \subset \Int Q(Z)$, $P \circ f(R_i) \subset \Int P(Z)$
 and $Q \circ f(R_i) = Q(Z)$ for any $i \in I$,
 \item the family $(\Int P(R_i))_{i \in I}$ covers $P(Z)$
 and its Lebesgue number is greater than $\theta \cdot \diam Q(Z)$,
 {\it i.e.}, for any subset $K$ of $P(Z)$ with
 $\diam K \leq \theta \diam Q(Z)$
 there exists $i \in I$ such that $K \subset \Int P(R_i)$.
\end{enumerate}
We call $Z$ and $\theta$
 the {\it superposition domain}
 and the {\it cone width} of the blending machine $\mc{B}$ respectively.
For a blending machine $\mc{B}=(R_i)_{i \in I}$ for a map $f$,
 we define the {\it forward invariant set} $\Lambda^s(\mc{B})$ by
$\Lambda^s(\mc{B})=\Lambda^s(f,\bigcup_{i \in I}R_i)$.
\begin{rmk}
\label{rmk:blender persistence}
A blending machine is $C^1$-persistent in the following sense:
If an open neighborhood $U_i$ of $R_i$
 and a compact subset $K_i$ of $\Int R_i$ are given for each $i \in I$
 then there exists a $C^1$-neighborhood $\cU$ of $f$
 such that any $g \in \cU$ admits a blending machine
 $\mc{B}'=(R'_i)_{i \in I}$ with respect to $(P,Q)$
 such that the superposition domain of $\mc{B}'$ is $Z$,
 the cone width is $\theta$,
 $R'_i \subset U_i$, and $K_i \subset \Int R'_i$ for any $i \in I$.
This fact follows from Lemma \ref{lemma:rectangle stability}
 and the openness of the conditions that defines a blending machine.
\end{rmk}
\begin{expl}
\label{expl:blender horseshoe}
This is essentially same as the affine model generating a blender
 explained in \cite[Section 6.2.1]{BDV}.
Define affine maps $f_1,f_2: \RR^3 \ra \RR^3$ by
\begin{alignat*}{4}
f_1(x_1,x_2,x_3)
 & =\left(\frac{1}{9}x_1-\frac{2}{3},\frac{7}{6}\left(x_2+1\right),
 8(x_3+1) \right), \\
f_2(x_1,x_2,x_3)
 & = \left(\frac{1}{9}x_1+\frac{2}{3},\frac{7}{6}\left(x_2-1\right),
 8(x_3-1) \right).
\end{alignat*}
Put
\begin{alignat*}{4}
K_1 & = [-3,3] \times [-13,11]  \times
 \left[-\frac{3}{2},-\frac{1}{2}\right], & \hsp
K_2 & = [-3,3] \times [-11,13]  \times
 \left[\frac{1}{2},\frac{3}{2}\right],
 \\
 R_1 & = [-2,2] \times \left[-\frac{5}{2},\frac{1}{2}\right]
  \times \left[-\frac{5}{4},-\frac{3}{4}\right], & \hsp
 R_2 & = [-2,2] \times \left[-\frac{1}{2},\frac{5}{2}\right]
  \times \left[\frac{3}{4},\frac{5}{4}\right].
\end{alignat*}
Remark that $R_i \subset \Int K_i$.
By direct computation. we have
\begin{alignat*}{4}
 f_1(K_1) & =\left[-1,-\frac{1}{3}\right]
 \times[-14,14] \times [-4,4] & \hsp
 f_1(K_2) & =\left[\frac{1}{3},1\right]
 \times[-14,14] \times [-4,4] \\
 f_1(R_1) & = \left[-\frac{8}{9},-\frac{4}{9}\right]
  \times \left[-\frac{7}{4},\frac{7}{4}\right]
  \times [-2,2], & \hsp
 f_2(R_2) & = \left[\frac{4}{9},\frac{8}{9}\right]
 \times \left[-\frac{7}{4},\frac{7}{4}\right]
  \times [-2,2].
\end{alignat*}
Since both $(K_1,K_2)$ and $(f_1(K_1),f_2(K_2))$ are pairs
 of mutually disjoint compact sets diffeomorphic to $[-1,1]^3$,
 we can take a $C^\infty$ diffeomorphism $f:\RR^3 \ra \RR^3$ such that
 $f(x)=f_i(x)$ for $i=1,2$ and $x \in K_i$.
See Figure \ref{fig:blender2}.
\begin{figure}
\begin{center}
 \includegraphics[scale=0.6]{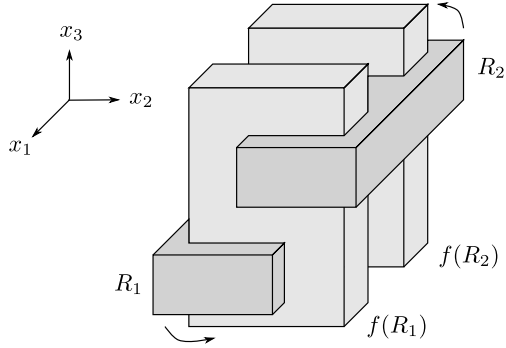} 
\end{center}
 \caption{An affine model generating a blender}
\label{fig:blender2}
\end{figure}
Put $U=M=\RR^3$, $Z=[-1,1] \times [-2,2]^2$,
 and define a $(2,1)$-splitting $(P,Q)$ of $M=\RR^3$
 by $P(x_1,x_2,x_3)=(x_1,x_2)$ and $Q(x_1,x_2,x_3)=x_3$.
Then, the family $(\Int P(R_i))_{i=1,2}$ covers $[-1,1] \times [-2,2]=P(Z)$
 and the Lebesgue number $\Delta$ of the cover is $1$.
The map $f_i$ satisfies the $\theta$-cone condition
 with respect to the splitting $(P,Q)$ for any $\theta>0$,
 for instance, $\theta=1/5<\Delta/\diam Q(Z)$.
The compact set $R_i$ is a $(P,Q \circ f)$-rectangle
 such that $Q(R_i) \subset (-2,2)=\Int Q(Z)$,
 $P(f(R_i)) \subset (-1,1) \times (-2,2) =\Int P(Z)$,
 $Q(f(R_i)) =[-2,2] =Q(Z)$ for each $i=1,2$.
Hence, $(R_i)_{i=1,2}$ is a blending machine for $f$
 with respect to the splitting $(P,Q)$
 such that the superposition domain is $Z$
 and the cone width $\theta$ for any $0<\theta<\Delta/\diam Q(Z)$.
The forward invariant set of the blending machine $\mc{B}$
 can be computed as
\begin{equation*}
 \Lambda^s(\mc{B})=\bigcup_{(x_1,x_2,x_3) \in \Lambda(R_1 \cup R_2,f)}
 [-2,2] \times \{(x_2,x_3)\} \subset W^s(\Lambda(K_1\cup K_2,f)).
\end{equation*}
Define a $(1,2)$-splitting $(P^h,Q^h)$ of $\RR^3$ by
 $P^h(x_1,x_2,x_3)=x_1$ and $Q^h(x_1,x_2,x_3)=(x_2,x_3)$.
It is easy to check that $f$, $(K_1,K_2)$, and $(P^h,Q^h)$
 satisfy the assumption of Proposition \ref{prop:transitive}.
This implies that $\Lambda(K_1 \cup K_2,f)$ is a
 locally maximal and topologically transitive hyperbolic set
 of unstable index two.
In fact, the hyperbolic splitting
 $T\RR^3|_{\Lambda(K_1 \cup K_2,f)}=E^s \oplus E^u$
 is given by $E^s(x)=\RR \oplus \{0\}$ and $E^u(x)=\{0\} \oplus \RR^2$
 for $x \in \Lambda(K_1 \cup K_2,f)$.
By Proposition \ref{prop:blender} below,
 for any $C^1$ map $\sigma:[-2,2] \ra (-1,1) \times (-2,2)$
 with $\|D\sigma\|<1/5$, the graph
 $\Gamma(\sigma)=\{(\sigma(t),t) \mid t \in [-2,2]\}$
 intersects with $\Lambda^s(\mc{B})$, and hence,
 it intersects with the stable set $W^s(\Lambda(K_1 \cup K_2, f),f)$.
This means that the hyperbolic set $\Lambda(K_1 \cup K_2,f)$ is a $cu$-blender
 of $uu$-index one.
\end{expl}
\begin{expl}[A `blender Cantor set' in $\RR^2$]
\label{expl:Cantor blender}
Define affine maps $f_1,f_2: \RR^2 \ra \RR^2$ by
\begin{alignat*}{4}
f_1(x,y) & =\left(\frac{7}{6}(x+1), 8(y+1) \right), &\hsp
f_2(x,y) & = \left(\frac{7}{6}(x-1), 8(y-1) \right)
\end{alignat*}
 and put
\begin{alignat*}{4}
 K_1 & =\left[-13,11 \right] \times \left[-\frac{3}{2},-\frac{1}{2}\right],
 & \hsp
 K_2 & =\left[-11,13 \right] \times  \left[\frac{1}{2},\frac{3}{2} \right],
\\
 R_1 & = \left[-\frac{5}{2},\frac{1}{2}\right]
 \times \left[-\frac{5}{4},-\frac{3}{4}\right], &\hsp
 R_2 & =  \left[-\frac{1}{2},\frac{5}{2}\right]
 \times \left[\frac{3}{4},\frac{5}{4}\right]
\end{alignat*}
See Figure \ref{fig:blender1}.
\begin{figure}
\begin{center}
 \includegraphics[scale=0.6]{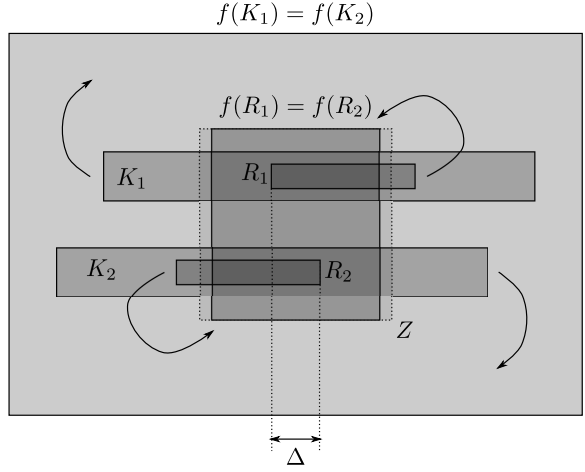} 
\end{center}
 \caption{The map $f$ generating a `blender Cantor set'}
\label{fig:blender1}
\end{figure}
Notice that $f_i$, $K_i$, and $R_i$ are the $(x_2,x_3)$-component of
 the map and the rectangles in the previous example.
We take a $C^1$ map $f:\RR^2 \ra \RR^2$ such that
 $f(x,y)=f_i(x,y)$ on a neighborhood of $K_i$.
Put $U=M=\RR^2$, $Z=[2,2]^2$ ,
 and define a $(1,1)$-splitting $(P,Q)$ of $\RR^2=M$
 by $P(x,y)=x$ and $Q(x,y)=y$.
In the same way as the previous example,
 one can check that $(R_i)_{i=1,2}$ is a blending machine
 for $f$ with respect to $(P,Q)$ such that
 the superposition domain is $Z$.
Since $K_1 \cup K_2 \subset (-14,14) \times (-4,4) = \Int f_i(K_i)$
 and $f_i$ is uniformly expanding on $K_i$ for $i=1,2$,
 the maximal forward invariant set $\Lambda^s(K_1 \cup K_2,f)$
 of $f$ in $K_1 \cup K_2$ is a $C^\infty$-regular Cantor set.
Since $R_i \subset \Int K_i$,
 the Cantor set $\Lambda^s(K_1 \cup K_2,f)$ contains $\Lambda^s(\mc{B})$.
Same as the previous example, for any $C^1$ function
 $\sigma:[-2,2] \ra (-2,2)$ with $|\sigma'|<1/5$, the graph
$\Gamma(\sigma)=\{(\sigma(t),t) \mid t \in [-2,2]\}$
intersects with $\Lambda^s(\mc{B})$, and hence,
it intersects with the Cantor set $\Lambda^s(K_1 \cup K_2,f)$.
\end{expl}
A blending machine induces a $cu$-blender in the following sense.
\begin{prop}
\label{prop:blender}
Let $(P,Q)$ be an $(l,m)$-splitting on an open manifold $M$
 and $f:U \ra M$ a $C^1$ embedding from an open subset $U$ of $M$
 into $M$.
Suppose that $f$ admits a blending machine $\mc{B}=(R_i)_{i \in I}$
 such that the superposition domain
 and the cone width are $Z$ and $\theta$ with $0<\theta<1$.
Let $\mc{F}$ be the set of $C^1$ maps $\sigma:Q(Z) \ra \Int P(Z)$ with
  $\|D\sigma\|<\theta$.
Define the `graph' $\Gamma(\sigma)$ of $\sigma \in \mc{F}$ by
\begin{equation*}
\Gamma(\sigma)=\{x \in Z \mid P(x)=\sigma(Q(x))\}.
\end{equation*}
Then, $\Gamma(\sigma)$ intersects with $\Lambda^s(\mc{B})$
 for any $\sigma \in \mc{F}$.

Moreover, if there exists a compact subset $K$ of $U$
 such that $R_i \subset \Int K$ for any $i \in I$
 and $\Lambda(K,f)$ is a hyperbolic invariant set of $f$
 with $\Lambda(K,f) \subset \Int K$,
 then $\Lambda(K,f)$ is a $cu$-blender of $uu$-index $m$
 whose superposition region is
 the set of $m$-dimensional $C^1$ embedded disks containing $\Gamma(\sigma)$
 in its interior for some $\sigma \in \mc{F}$.
\end{prop}
\begin{proof}
We prove the proposition by a `graph transformation argument'
 introduced in \cite{BKR14}.
Let $\Delta$ be the Lebesgue number of the open cover
 $(\Int P(R_i))_{i \in I}$ of $P(Z)$.
Take $\sigma_0 \in \mc{F}$.
Then, $\Gamma(\sigma_0) \subset Z$ and
 $\diam P(\Gamma(\sigma_0))<\theta \diam Q(Z)<\Delta$.
Hence, there exists $i_0 \in I$ such that
 $P(\Gamma(\sigma_0)) \subset \Int P(R_{i_0})$.
Since $Q(R_{i_0}) \subset \Int Q(Z)=\Int Q(\Gamma(\sigma_0))$,
we can show that the `graph' $\Gamma(\sigma_0)$ intersects with
 $D^-=\{y \in R_{i_0} \mid Q \circ f(x)=t_0\}$
 at a unique point for any $t_0 \in Q(Z)$
 by the same argument as Lemma \ref{lemma:rectangle composition}.
This implies that there exists a $C^1$ map $\sigma_1:Q(Z) \ra P(Z)$
 such that $f(\Gamma(\sigma_0) \cap R_{i_0})=\Gamma(\sigma_1)$.
Since $P(f(R_{i_0})) \subset \Int P(Z)$
 and $f$ satisfies the $\theta$-cone condition with respect to $(P,Q)$,
 one can see that the map $\sigma_1$ is an element of $\mc{F}$.
Iterating this process, we obtain sequences $(i_n)_{n \geq 0}$ in $I$
 and $(\sigma_n)_{n \geq 0}$ in $\mc{F}$ such that
 $f(\Gamma(\sigma_n) \cap R_{i_n})=\Gamma(\sigma_{n+1})$ for any $n \geq 0$.
Then,
\begin{equation*}
 \emptyset \neq \bigcap_{n \geq 0}f^{-n}(\Gamma(\sigma_n))
 \subset \Gamma(\sigma_0) \cap \bigcap_{n \geq 0} f^{-n}(R_{i_n})
 =\Gamma(\sigma_0) \cap \Lambda^s(\mc{B}).
\end{equation*}

Suppose that there exists a compact subset $K$ of $M$
 such that $R_i \subset \Int K$ for any $i \in I$
 and $\Lambda(K,f) \subset \Int K$.
Then, $\Lambda^s(\mc{B}) \subset \Lambda^s(K,f) \subset W^s(\Lambda(K,f))$
 and hence, $W^s(\Lambda(K,f))$ intersects with $\Gamma(\sigma)$
 for any $\sigma \in \mc{F}$.
By the persistence of blending machine,
 the same holds for any map $g$ which is $C^1$-sufficiently close to $f$.
Therefore, $\Lambda(K,f)$ is a $cu$-blender of $uu$-index $m$
 when $\Lambda(K,f)$ is a hyperbolic set of $f$.
\end{proof}

The following criterion to
 $C^1$-stable intersection of the forward invariant sets
 of blending machines for two maps
 is a keystone to prove Theorems A and B.
\begin{thm}
\label{thm:intersection}
Let $l,m$ be positive integers
 and $M_1, M_2$ be $(l+m)$-dimensional manifolds
 with $(l,m)$-splitting $(P_1,Q_2)$, $(m,l)$-splitting $(P_2,Q_2)$
 respectively.
For $\tau=1,2$,
 let $U_\tau$ be an open subset of $M_\tau$,
 $f_\tau:U_\tau \ra M_\tau$ a $C^1$ map,
 and $\mc{B}_\tau=(R_{\tau,i})_{i \in I_\tau}$ 
 is a blending machine for $f_\tau$ with respect to $(P_\tau,Q_\tau)$
 whose superposition domain is $Z_\tau$.
Suppose that there exist positive real numbers
 $\theta$, $\lambda_\tau,\mu_\tau$ ($\tau=1,2,\sh$),
 a compact subset $R_\sh$ of $U_2$,
 and a $C^1$ embedding $f_\sh:R_\sh \ra U_1$
 which satisfy the following properties:
\begin{enumerate}
\item
 \label{item:intersection 4}
 $\lambda_1\mu_2>1$, $\lambda_2\mu_1>1$.
\item
 \label{item:intersection 1}
 $f_\tau$ satisfies
 the $(\theta,\lambda_\tau,\mu_\tau)$-cone condition on $R_{\tau,i}$
 with respect to $(P_\tau, Q_\tau)$ for any $\tau=1,2$ and $i \in I_\tau$.
\item
 \label{item:intersection 3}
 $R_\sh$ is a $(Q_2, Q_1 \circ  f_\sh)$-rectangle with
\begin{alignat*}{4}
 P_2(R_\sh) & \subset \Int P_2(Z_2), &\hsp
 Q_2(R_\sh) & =Q_2(Z_2), \\
 P_1(f_\sh(R_\sh)) & \subset \Int P_1(Z_1), & \hsp
 Q_1(f_\sh(R_\sh)) & = \Int Q_1(Z_1).
\end{alignat*}
\item
 \label{item:intersection 2}
 $f_\sh$ satisfies
 the $(\theta,\lambda_\sh,\mu_\sh)$-cone condition on $R_\sh$
 with respect to $(Q_2,P_2)$ and $(P_1,Q_1)$.
\item
 \label{item:intersection 5}
 Let $\Delta_\tau$ be the Lebesgue number
 of the cover $(\Int P_\tau(R_{\tau,i}))_{i \in I_\tau}$ of $P_\tau(Z_\tau)$
 for $\tau=1,2$.
Then,
\begin{gather*}
 \theta \diam Q_1(Z_1)+\lambda_\sh^{-1}\diam Q_2(Z_2) < \Delta_1,\\
 \theta \diam Q_2(Z_2)+\mu_\sh^{-1}\diam Q_1(Z_1) < \Delta_2.
\end{gather*}
\end{enumerate}
Then, $\Lambda^s(\mc{B}_1)$
 intersects with $f_\sh(R_\sh \cap \Lambda^s(\mc{B}_2))$.

Moreover, for any compact neighborhoods $K_1$ and $K_2$
 of $\bigcup_{i \in I_1}R_{1,i}$ and $\bigcup_{j \in I_2}R_{2,j}$
 respectively,
 there exist neighborhoods $\cU_,\cU_2,\cU_\sh$
 of $f_1,f_2,f_\sh$
 with respect to the $C^1$ compact open topology
 such that
 $\Lambda^s(g_1,K_1) \cap g_\sh(\Lambda^s(g_2,K_2)) \neq \emptyset$
 for any $g_1 \in \cU_1$, $g_2 \in \cU_2$, and $g_\sh \in \cU_\sh$.
\end{thm}
\begin{proof}
We say that a pair $(R,h)$ of a compact subset of $U_2$ and
 a $C^1$ embedding $h:R \ra U_1$ satisfies the condition $(\sh)$
 if $R$ is a $(Q_2,Q_1 \circ h)$-rectangle
 with 
\begin{alignat*}{4}
 P_2(R) & \subset P_2(Z_2), & \hsp
 Q_2(R) & =Q_2(Z_2), \\ 
 P_1(h(R)) & \subset \Int P_1(Z_1), &\hsp
 Q_1(h(R)) & = \Int Q_1(Z_1),
\end{alignat*}
 and
 $h$ satisfies the $(\theta,\lambda_\sh,\mu_\sh)$-cone condition on $R$
 with respect to $(Q_2,P_2)$ and $(P_1,Q_1)$.
Notice that the pair $(R_\sh,f_\sh)$ satisfies the condition $(\sh)$.
We claim that if a pair $(R,h)$ satisfies the condition $(\sh)$
 then there exist $i \in I_1$ and $j \in I_2$ such that
 the pair $(f_2(h^{-1}(R_{1,i}) \cap R \cap R_{2,j}),
 f_1 \circ h \circ f_2^{-1})$ satisfies the condition $(\sh)$.

Before the proof of the claim,
 we see that the claim implies that
  $\Lambda^s(\mc{B}_1)$ intersects 
 with $f_\sh(R_\sh \cap \Lambda^s(\mc{B}_2))$.
When the claim holds,
 we can take sequences $(R_n,h_n)_{n \geq 0}$, 
 $(i_n)_{n \geq 0}$, and $(j_n)_{n \geq 0}$ of
 pairs satisfying the condition $(\sh)$
 and elements of $I_1$, $I_2$ respectively
 such that $(R_0,h_0)=(R_\sh,f_\sh)$ and
\begin{alignat*}{4}
 R_{n+1} & =
 f_2(h_n^{-1}(R_{1,i_n}) \cap R_n \cap R_{2,j_n}), & \hsp
 h_{n+1} & =f_1 \circ h_n \circ f_2^{-1}
\end{alignat*}
 for any $n \geq 0$.
Then, $h_n  = f_1^n \circ f_\sh \circ f_2^{-n}$ and
\begin{gather*}
f_\sh \circ f_2^{-n}(R_n)
 =\bigcap_{k=0}^{n-1} f_1^{-k}(R_{1,i_k})
  \cap f_\sh\left(R_\sh \cap \bigcap_{k=0}^{n-1} f_2^{-k}(R_{2,j_k})\right).
\end{gather*}
Since $Q_2(R_n)=Q_2(Z_2)$,
 $R_n$ is a non-empty compact subset of $U_2$.
This implies that the intersection
\begin{equation*}
 \bigcap_{n \geq 0} f_\sh \circ f_2^{-n}(R_n)
 =  \bigcap_{n \geq 0}f_1^{-n}(R_{1,i_n})
 \cap f_\sh\left(R_\sh \cap \bigcap_{n \geq 0}f_2^{-n}(R_{2,j_n})\right)
\end{equation*}
 is non-empty.
The right-hand side is
 a subset of $\Lambda(\mc{B}_1) \cap f_\sh(\Lambda(\mc{B}_2))$.
Therefore, $\Lambda(\mc{B}_1) \cap f_\sh(\Lambda(\mc{B}_2))$ is non-empty.

Let us prove the claim above.
Suppose that a pair $(R,h)$ satisfies the condition $(\sh)$.
Notice that $Q_2(R)=Q_2(Z_2)$, $Q_1(h(R))=Q_1(Z_1)$,
 and $h$ satisfies the $(\theta,\lambda_\sh,\mu_\sh)$-condition
 with respect to the splittings $(Q_2,P_2)$ and $(P_1,Q_1)$.
By Corollary \ref{cor:rectangle size},
\begin{align*}
 \diam P_2(R) & \leq \theta\diam Q_2(Z_2) +\mu_\sh^{-1} Q_1(Z_1)<\Delta_2,\\
 \diam P_1(h(R)) & \leq \theta\diam Q_1(Z_1)
 +\lambda_\sh^{-1} Q_2(Z_2)<\Delta_1.
\end{align*}
Since $\Delta_\tau$ is the Lebesgue number of
 the open cover $(\Int P_\tau(R_{\tau,i}))_{i \in I_\tau}$
 of $P_\tau(Z_\tau)$,
 there exits $i_\sh \in I_1$ and $j_\sh \in I_2$
 such that $P_1(h(R)) \subset P_1(\Int R_{1,i_\sh})$
 and $P_2(R) \subset P_2(\Int R_{2,j_\sh})$.
By the definition of a blending machine, we also have
\begin{alignat*}{4}
 Q_1(R_{1,i_\sh}) &\subset \Int Q_1(Z_1)= \Int Q_1(h(R)),  & \hsp
 Q_2(R_{2,j_\sh}) &\subset \Int Q_2(Z_2)= \Int Q_2(R).
\end{alignat*}
We apply Lemma \ref{lemma:rectangle composition}
 to a $(Q_2,P_2 \circ f_2^{-1})$-rectangle $f_2(R_{2,j_\sh})$,
 a $(Q_2,Q_1 \circ h)$-rectangle $R$,
 a $(P_1,Q_1 \circ f_1)$-rectangle $R_{1,i_\sh}$,
 and maps $f_2^{-1}$, $h$, and $f_1$.
Then, a compact subset 
\begin{equation*}
 R'=f_2(h^{-1}(R_{1,i_\sh})) \cap f_2(R) \cap f_2(R_{2,j_\sh})
 =f_2(h^{-1}(R_{1,i_\sh}) \cap R \cap R_{2,j_\sh}) 
\end{equation*}
 of $M_2$ is a $(f_1 \circ h \circ f_2^{-1})$-rectangle with
\begin{gather*}
 Q_2(R') =Q_2(f_2(R_{2,j_\sh}))=Q_2(Z_2), \\
 Q_1(f_1 \circ h \circ f_2^{-1}(R'))
  =Q_1(f_1(R_{1,i_\sh}))= Q_1(Z_1).
\end{gather*}
By the cone condition on $f_1$, $f_2$, and $f_\sh$,
 the composition $f_1 \circ h \circ f_2^{-1}$ satisfies
 the $(\theta,\lambda_1\lambda_\sh\mu_2,\mu_1\mu_\sh\lambda_2)$-cone
 condition with respect to the splittings $(Q_2,P_2)$ and $(P_1,Q_1)$
 on $R'$.
Since $\lambda_1\mu_2>1$ and $\mu_1\lambda_2>1$,
 the map $f_1 \circ h \circ f_2^{-1}$ satisfies
 the $(\theta,\lambda_\sh,\mu_\sh)$-cone
 condition with respect to the splittings $(Q_2,P_2)$ and $(P_1,Q_1)$
 on $R'$.
Therefore, the pair $(R',f_1 \circ h \circ f_2^{-1})$ satisfies
 the condition $(\sh)$.
The proof of the claim is finished,
 and hence, $\Lambda^s(\mc{B}_1)$ intersects with
 $f_\sh(R_\sh \cap \Lambda^s(\mc{B}_2))$.
\begin{figure}
\begin{center}
 \includegraphics[scale=0.6]{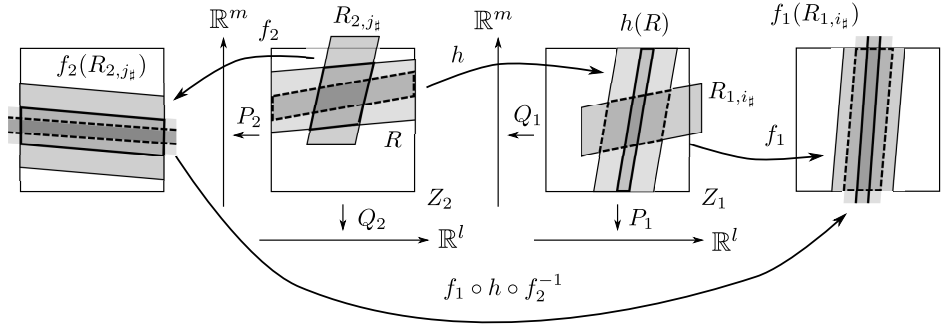} 
\end{center}
 \caption{Claim in Proof of Theorem \ref{thm:intersection}}
\end{figure}

Finally, we show the stability of intersection
 of the maximal forward invariant sets.
For $\tau=1,2$, let $K_\tau$ be
 a compact neighborhood of $\bigcup_{i \in I_\tau}R_{\tau,i}$.
By Lemma \ref{lemma:rectangle stability}
 and the persistence of the cone condition,
 there exist $C^1$-neighborhoods
 $\cU_1$, $\cU_2$, and $\cU_\sh$ of $f_1$, $f_2$, and $f_3$ respectively
 such that any triple $(g_1,g_2,g_3)$ 
 in $\cU_1 \times \cU_2 \times \cU_\sh$
 admits a blending machine $\mc{B}'_\tau=(R'_{\tau,i})_{i \in I_\tau}$
 and a compact set $R'_\sh$ which satisfy
 the assumption of the theorem
 and $\bigcup_{i \in I_\tau}R'_{\tau,i} \subset \Int K_\tau$
 for $\tau=1,2$.
This implies that
\begin{equation*}
 \Lambda^s(K_1,g_1) \cap g_\sh(\Lambda^s(K_2,g_2))
 \supset \Lambda^s(\mc{B}'_1)
 \cap g_\sh(R'_\sh \cap \Lambda^s(\mc{B}'_2))  \neq \emptyset.
\end{equation*}
 for any $(g_1,g_2,g_\sh) \in \cU_1 \times \cU_2 \times \cU_\sh$.
\end{proof}

\section{Stable intersection of Cantor sets}
\label{sec:Cantor}
In this section, we prove Theorem A stated in the introduction.
\begin{thmA}
For any integers $l,m \geq 1$ and any real number $\delta>0$,
 there exists a pair $(\Lambda_1,\Lambda_2)$
 of $C^\infty$-regular Cantor sets in $\RR^{l+m}$
 which has $C^1$-stable intersection
 and satisfy that $\udim_B(\Lambda_1)<l+\delta$,
 $\udim_B(\Lambda_2)<m+\delta$.
\end{thmA}

Before the proof,
 we give the precise definitions of stable intersection
 of a pair of regular Cantor sets which belong
 to the class of Cantor sets given
 in Example \ref{expl:regular Cantor}\footnote{
We can define stability of intersections of
 a pair of general regular Cantor sets in a similar way.
 But, we omit it here in order to simplify notations
 since the Cantor sets we construct in the proof of Theorem A belong
 to the class of Cantors sets in Example \ref{expl:regular Cantor}.}.
Fix a finite set $I$ with $\# I \geq 2$
 and a family $\mc{K}=(K_i)_{i \in I}$ of mutually disjoint
 non-empty compact subset of $\RR^m$.
Put $K=\bigcup_{i \in I} K_i$.
Let $\mc{C}^r(\mc{K})$ be the set of $C^r$ self-maps $f$ on $\RR^m$
 such that $K \subset \Int f(K_i)$ and
 the restriction of $f$ to $K_i$ is uniformly expanding for any $i \in I$.
This set is an open subset of the space
 $C^r(\RR^m,\RR^m)$ of $C^r$ self-maps on $\RR^m$.
As in Example \ref{expl:regular Cantor},
 the maximal forward invariant set $\Lambda^s(K,f)$ is
 a $C^r$-regular Cantor set for $f \in \mc{C}^r(\mc{K})$.
Fix $1 \leq r'\leq r \leq \infty$, finite sets $I,J$,
 and families $\mc{K}_1=(K_{1,i})_{i \in I}$,
 $\mc{K}_2=(K_{2,j})_{j \in J}$ such that each $\mc{K}_\tau$
 is a family of mutually disjoint non-empty compact subsets of $\RR^m$.
Put $K_1=\bigcup_{i \in I}K_{1,i}$ and $K_2=\bigcup_{j \in J}K_{2,j}$.
For $f \in \mc{C}^r(\mc{K}_1)$ and $g \in \mc{C}^r(\mc{K}_1)$,
 we say that the pair $(\Lambda^s(K_1,f),\Lambda^s(K_2,g))$
 of $C^r$-regular Cantor sets has {\it $C^{r'}$-stable intersection}
 if there exist $C^{r'}$-neighborhoods
 $\mc{U}$ of $f$ in $\mc{C}^{r'}(\mc{K}_1)$ and
 $\mc{V}$ of $g$ in $\mc{C}^{r'}(\mc{K}_2)$
 such that $\Lambda^s(K_1,f')$ intersects with $\Lambda^s(K_2,g')$
 for any $f' \in \mc{U}$ and $g' \in \mc{V}$.

For $m \geq 1$, $p \in \RR^m$, and $r>0$,
 let $B^m(p,r)$ be the $m$-dimensional closed $r$-ball 
 $\{x \in \RR^m \mid \|x-p\| \leq r\}$ with respect to the box norm.
When $p$ is the origin of $\RR^m$, we write just $B^m(r)$ for it.
Remark that $B^l(p,r)$ is diffeomorphic to $[0,1]^l$,
 $B^l(r)=[-r,r]^l$,
 and $B^l(p,r) \times B^m(q,r)=B^{l+m}((p,q),r)$
 for $p \in \RR^l$, $q \in \RR^m$, and $r>0$.
The following pasting lemma is classical
 and easy to prove:
\begin{lemma}
\label{lemma:paste} 
Let $(R_i)_{i \in I}$ a family of mutually disjoint
 compact subsets of $\RR^m$
 which are $C^r$ diffeomorphic to $[0,1]^m$
 and $(f_i)_{i \in I}$ a $C^r$ orientation preserving
 diffeomorphisms of $\RR^m$
 with $m \geq 2$, $1 \leq r \leq \infty$, and a finite set $I$.
Then, there exists a $C^r$ self-map $f$ on $\RR^m$
 such that $f=f_i$ on a neighborhood of $R_i$ for each $i \in I$
 and the support of $f$, {\it i.e.},
 the closure of ${\{x \in \RR^m \mid f(x) \neq x\}}$, is compact.
Moreover,
 if the family $(f_i(R_i))_{i \in I}$ consists
 of mutually disjoint compact sets,
 then we can take the map $f$ above
 so that $f$ is a $C^r$ diffeomorphism of $\RR^m$.
\end{lemma}
Let us prove Theorem A.
It is done by constructing two diffeomorphisms
 $f_1,f_2$ with blending machines and the `connecting map' $f_\sh$
 which satisfy the assumption of Theorem \ref{thm:intersection}.
\begin{proof}
[Proof of Theorem A]
Fix $l,m \geq 1$ and $\delta>0$.
Put $(l_1,m_1)=(l,m)$ and $(l_2,m_2)=(m,l)$.
For $\tau=1,2$,
 let $(P_\tau,Q_\tau):\RR^{l+m} \ra \RR^{l_\tau} \times \RR^{m_\tau}$
 be the natural $(l_\tau,m_\tau)$-splitting.
Take a finite set $I_\tau$
 and a family $(p^-_{\tau,i})_{i \in I_\tau}$ of points
 in $\Int B^{l_\tau}(2)$
 such that the family $(\Int B^{l_\tau}(p^-_{\tau,i},1))_{i \in I_\tau}$
 is an open cover of $B^{l_\tau}(2)$ for each $\tau=1,2$.
We denote the cardinality of $I_\tau$ by $\# I_\tau$
 and the Lebesgue number of the open cover by $\Delta_\tau$.
Take a family $(p^+_{\tau,i})_{i \in I_\tau}$ of
 mutually distinct points in $\Int B^{m_\tau}(1)$
 and choose a real number $\mu_\tau> \max\{2,(\# I)^{1/\delta}\}$
 such that $(B^{m_\tau}(p^+_{\tau,i},1/\mu_\tau))_{i \in I_\tau}$
 is a family of mutually disjoint subsets of $\Int B^{m_\tau}(1)$.
For $i \in I_\tau$, we put
\begin{alignat*}{4}
 R_{\tau,i} & = B^{l_\tau}(p^-_{\tau,i},1)
  \times B^{m_\tau}(p^+_{\tau,i},1/\mu_\tau), & \hsp
 K_{\tau,i} & = B^{l_\tau}(p^-_{\tau,i},6)
 \times B^{m_\tau}(p^+_{\tau,i},2/\mu_\tau).
\end{alignat*}
 Since the compact sets in the family $(K_{\tau,i})_{i \in \tau}$ 
 are diffeomorphic to $[-1,1]^{l+m}$
 and mutually disjoint for each $\tau=1,2,$
 Lemma \ref{lemma:paste} implies that
 there exist $C^\infty$ self-maps $f_1$ and $f_2$ on $\RR^{l+m}$ such that
\begin{equation*}
 f_\tau(x,y)=\left(\frac{3}{2}(x-p^-_{\tau,i}),
 \mu_\tau(y-p^+_{\tau,i})\right)
\end{equation*}
 for $\tau=1,2$, $i \in I_\tau$, and $(x,y) \in K_{\tau,i}$.
Put $K_\tau=\bigcup_{i \in I_\tau} K_{\tau,i}$
 and $Z_\tau=B^{l_\tau}(2) \times B^{m_\tau}(1)$.
Notice that
 the family $(P_\tau(R_{\tau,i}))_{i \in I_\tau}$ covers
 $P_\tau(Z_\tau)=B^{l_\tau}(2)$.
The compact set $R_{\tau,i}$
 is a $(P_\tau,Q_\tau \circ f_\tau))$-rectangle
 and $f_\tau(R_{\tau,i})=B^{l_\tau}(3/2) \times B^{m_\tau}(1)$.
In particular, $P_\tau(f_\tau(R_{\tau,i})) \subset P_\tau(Z_\tau)$
 and $Q_\tau(f_\tau(R_{\tau,i}))=Q_\tau(Z_\tau)$.
Take $\theta>0$ such that
 $\theta \diam Q_\tau(Z_\tau)<\Delta_\tau$ for $\tau=1,2$.
Then, the map $f_\tau$ satisfies
 the $(\theta,1,2)$-cone condition on $K_\tau$
 for each $i \in I_\tau$.
Therefore, the family $\mc{B}_\tau=(R_{\tau,i})_{i \in I_\tau}$
 is a blending machine for $f_\tau$ with respect to $(P_\tau,Q_\tau)$
 such that the superposition domain and the cone width are
 $Z_\tau$ and $\theta$.

Choose large numbers $\lambda_\sh,\mu_\sh>1$ so that
\begin{gather*}
 \theta\diam Q_1(Z_1)+\lambda_\sh^{-1}\diam Q_2(Z_2) < \Delta_1, \\
 \theta\diam Q_2(Z_2)+\mu_\sh^{-1}\diam Q_1(Z_1) < \Delta_2.
\end{gather*}
Put $R_\sh=B^m(\mu_\sh^{-1}) \times B^l(1)$
 and define an affine map $f_\sh:\RR^{l+m} \ra \RR^{l+m}$ by
\begin{equation*}
 f_\sh(x,y)=(\lambda_\sh^{-1}y, \mu_\sh x)
\end{equation*}
 for $x \in \RR^m$ and $y \in \RR^l$.
Then, 
\begin{alignat*}{4}
P_2(R_\sh) & \subset \Int B^m(2) = \Int P_2(Z_2), &\hsp
Q_2(R_\sh) & = B^l(1) =Q_2(Z_2),\\
P_1(f_\sh(R_\sh)) & \subset \Int B^l(2) = \Int P_1(Z_1), &\hsp
Q_1(f_\sh(R_\sh)) & = B^m(1) = Q_1(Z_1).
\end{alignat*}
It is easy to check that
 $R_\sh$ is a $(Q_2,Q_1 \circ f_\sh)$-rectangle
 and
 $f_\sh$ satisfies the $(\theta,\lambda_\sh,\mu_\sh)$-cone condition
 with respect to $(Q_2,P_2)$ and $(P_1,Q_1)$.
Therefore, the maps $f_1$, $f_2$, $f_\sh$,
 the families $\mc{B}_1=(R_{1,i})_{i \in I_1}$,
 $\mc{B}_2=(R_{2,j})_{j \in I_2}$,
 and the compact set $R_\sh$ satisfy the assumption of
 Theorem \ref{thm:intersection}.
Since $R_{\tau,i} \subset \Int K_{\tau,i}$,
 for any $\tau=1,2$ and $i \in I_\tau$,
 Theorem \ref{thm:intersection} implies
 that $\Lambda^s(K_1,g_1)$ intersects with
 $f_\sh(\Lambda^s(K_2,g_2))$
 for any $C^1$ maps $g_1,g_2$ which are $C^1$-close to $f_1, f_2$ respectively.
The map $f_\tau$ is uniformly expanding on each $K_{\tau,i}$
 and $K_\tau \subset \Int f_\tau(K_{\tau,i})
 = \Int (B^{l_\tau}(9) \times B^{m_\tau}(2))$
 for each $\tau=1,2$ and $i \in I_\tau$.
Hence, $\Lambda^s(K_\tau,g_\tau)$ is a $C^r$-regular Cantor set
 for any $g_\tau$ sufficiently $C^1$ close to $f_\tau$.

Put $f_2^T=f_\sh \circ f_2 \circ f_\sh^{-1}$.
Then,
\begin{equation*}
 f_2^T(x,y)
 =\left(\frac{3}{2}(x-\lambda_\sh^{-1}p^+_{2,j}),
 \mu_2(y-\mu_\sh p^-_{2,j})\right)
\end{equation*}
 for $(x,y) \in f_\sh(K_{2,j})$,
 and hence, the map $f_2^T$ is expanding on each $f_\sh(K_{2,j})$.
Hence, the maximal forward invariant set
 $\Lambda^s(f_\sh(K_2),f_2^T)=f_\sh(\Lambda^s(K_2,f_2))$
 is a $C^\infty$-regular Cantor set.
For any $C^1$ map $g^T_2$ which is $C^1$-close to $f^T_2$,
 the map $f_\sh^{-1} \circ g^T_2, \circ f_\sh$ is
 $C^1$-close to $f_2$ and
\begin{equation*}
 \Lambda^s(f_\sh(K_2),g^T_2)
 =f_\sh(\Lambda^s(K_2,f_\sh^{-1} \circ g^T_2, \circ f_\sh)).
\end{equation*}
This implies that
 $\Lambda^s(K_1,g_1)$ intersects with $\Lambda^s(f_\sh(K_2),g^T_2)$
 for any maps $g_1$ and $g^T_2$ which are $C^1$-close to
 $f_1$ and $f^T_2$.
Therefore, the $C^\infty$-regular Cantor sets
 $\Lambda^s(K_1,f_1)$ and $\Lambda^s(f_\sh(K_2),f_2^T)$
 have $C^1$-stable intersection.
\begin{figure}
\begin{center}
 \includegraphics[scale=0.6]{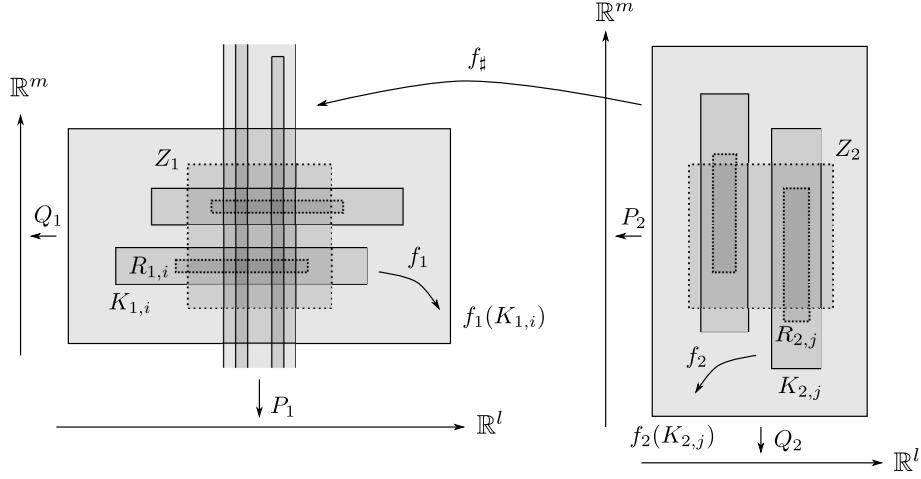} 
\end{center}
 \caption{Proof of Theorem A}
\label{fig:Cantor}
\end{figure}

Finally, we compute the upper box dimensions of $\Lambda^s(K_1,f_1)$
 and $\Lambda^s(f_\sh(K_2),f^T_2)=f_\sh(\Lambda^s(K_2,f_2))$.
Let $f_\tau^+:\RR^{m_\tau} \ra \RR^{m_\tau}$ be a $C^\infty$ map such that
 $f_\tau^+(y)=\mu_\tau(y-p^+_{\tau,i})$
 on $B^{m_\tau}(p^+_{\tau,i},2/\mu_\tau)$.
Put $K_\tau^+=\bigcup_{i \in I_\tau}B^m(p^+_{\tau},2/\mu_\tau)$.
Since $f_\tau(B^{l_\tau}(p^-_{\tau},6) \times \{y\})
 \subset B^{l_\tau}(9) \times \{f_\tau^+(y)\}$
 for any $i \in I_\tau$ and $y \in B^{m_\tau}(p^+_{\tau,i},2/\mu_\tau)$,
 the set $\Lambda^s(K_\tau,f_\tau)$
 is a subset of $B^{l_\tau}(6) \times \Lambda^s(K_\tau^+,f_\tau^+)$.
For any $N \geq 1$,
 $\bigcap_{n=0}^{N-1}(f_\tau^+)^{-n}(K_\tau^+)$ is
 the union of mutually disjoint $(\# I_\tau)^N$
 balls of radius $2/\mu_\tau^N$.
It follows that
 the upper box dimension of
 $\Lambda^s(K_\tau^+,f_\tau^+)$ is
 $-\log \# I/\log(1/\mu_\tau)=\log\# I/\log \mu_\tau$
 by the standard computation.
By the product inequality of the upper box dimension,
\begin{equation*}
 \udim_B \Lambda^s(K_\tau,f_\tau) \leq \udim_B B^{l_\tau}(9)
 + \udim \Lambda^s(K_\tau^+,f^+_\tau)
 < l_\tau +\frac{\log \# I_\tau}{\log \mu_\tau}<l_\tau+\delta.
\end{equation*}
Hence, $\udim_B \Lambda^s(K_1,f_1)<l+\delta$.
It is known that any diffeomorphism of $\RR^{l+m}$
 preserves the upper box dimension of a compact set.
Hence,
 $\udim_B \Lambda^s(f_\sh(K_2),f_2^T) =\udim_B f_\sh(\Lambda^s(K_2,f_2))
 <l_2+\delta=m+\delta$.
\end{proof}

\section{Robust homoclinic tangency of the largest codimension}
\label{sec:robust}
In this section, we prove
\begin{thm}
\label{thm:ThmB'}
For any $m \geq 2$,
 there exists a $C^\infty$ diffeomorphism of $\RR^{2m}$
 whose support is compact and
 which exhibits $C^2$-robust homoclinic tangency of codimension $m$.
\end{thm}
We can obtain a diffeomorphism of $\RR^{2m+1}$
 exhibiting $C^2$-robust homoclinic tangency of codimension $m$
 by taking a product with a one-dimensional strong contraction.
To prove Theorem B for a general manifold of dimension at least four,
 it is sufficient to embed such a dynamics into the given manifold.

\subsection{Lifts of diffeomorphisms}
\label{sec:lift}
Fix integers $l,m \geq 1$.
We denote the set of real $(m,l)$-matrices by $\Mat(m,l)$.
Let $M$ be an $(l+m)$-dimensional manifold.
By $\Gr_x(M,l)$, we denote the set of $l$-dimensional linear subspaces
 of the tangent space $T_x M$ at $x \in M$.
The Grassmannian bundle $\Gr(M,l)=\bigcup_{x \in M} \Gr_x(M,l)$
 is a $C^\infty$ fiber bundle 
 with $lm$-dimensional fibers $\Gr_x(M,l)$.
We denote the natural projection from $\Gr(M,l)$ to $M$ by $\pi$.
The smooth structure of $\Gr(M,l)$ is given as follows:
For $x \in M$ and $\xi \in \Gr_x(M,l)$,
 take a smooth local coordinate $(U,\vphi)$ of $M$ at $x$
 such that $D\vphi(\xi)=\RR^l \oplus \{0\}$.
Put
\begin{equation*}
 V_\vphi
 = \{\eta \in \Gr(M,l) \mid \pi(\eta) \in U,
 D\vphi(\eta) \text{ is transverse to } \{0\} \oplus \RR^m\}.
\end{equation*}
For each $\eta \in V_\vphi$,
 there exists $\Pi_\vphi(\eta) \in \Mat(m,l)$
 such that $D\vphi(\eta)=\{(v,\Pi_\vphi(\eta)(v)) \mid v \in \RR^{l}\}$.
The map $(\pi,\Pi_\vphi):\eta \mapsto (\pi(\eta),\Pi_\vphi(\eta))$
 is a bijection from $V_\vphi$ to $\pi(U) \times \Mat(m,l)$ and
 the family $\{(V_\vphi, (\pi,\Pi_\vphi))\}$ gives
 the smooth structure of $\Gr(M,l)$ as a fiber bundle.

For an open subset $U$ of $M$ and an $C^1$ embedding $f:U \ra M$,
 we define the {\it lift} $\wh{f}:\Gr(U,l) \ra \Gr(M,l)$ of $f$ by
 $\wh{f}(\xi)=Df(\xi)=\{Df v \mid v \in \xi\}$ for $\xi \in \Gr(U,l)$.
Remark that the lift $\wh{f}$ is of class $C^{k-1}$ if $f$ is of class $C^k$.
\begin{prop}
\label{prop:lift W-s}
Let $M$ be an $(l+m)$-dimensional  manifold
 with an $(l,m)$-splitting $(P^h,Q^h)$,
 $K$ a compact subset of $M$
 and $f:K \ra M$ a $C^2$-embedding.
Suppose that $\Lambda(K,f) \subset \Int K$ and
 $f$ satisfies the $(\theta,1,1)$-cone condition on $K$
 with respect to the splitting $(P^h,Q^h)$
 (by Remark \ref{rmk:cone hyp}, $\Lambda$ is a hyperbolic set
 of stable index $l$).
Put 
\begin{equation*}
 \wh{K}=\{\xi \in \pi^{-1}(K) \mid \xi \subset C(\pi(\xi),\theta,P^h,Q^h)\}.
\end{equation*}
 and let $\wh{f}:\pi^{-1}(K) \ra \Gr(M,l)$ be the lift of $f$
 to the Grassmannian bundle.
Then,
\begin{equation*}
 \Lambda^s(\wh{K},\wh{f})
 \subset \{T_q W^s(p) \mid p \in \Lambda(K,f),
 q \in W^s(p,f) \cap \Lambda^s(K,f)\}.
\end{equation*}
\end{prop}
\begin{proof}
Take $\xi \in \Lambda^s(\wh{K},\wh{f})$ and put $q=\pi(\xi)$.
Then, $f^n(q)=\pi(\wh{f}^n(\xi))$
 is contained in $\pi(\wh{K})=K$ for any $n \geq 0$,
 and hence, $q$ is an element of $\Lambda^s(K,f)$.
By compactness of $K$,
 $q$ is contained in the stable set $W^s(\Lambda(K,f),f)$.
Since $\Lambda(K,f)$ is a compact locally maximal hyperbolic invariant set
 of $f$, the stable set $W^s(\Lambda(K,f),f)$
 coincides with the union of the stable manifolds $W^s(p,f)$ of
 points $p$ in $\Lambda(K,f)$.
Hence, there exists $p \in \Lambda(K,f)$
 such that $q \in W^s(p,f) \cap \Lambda^s(K,f)$.
The proof will finish once we show that $\xi \subset T_q W^s(p,f)$
 since $\dim \xi=\dim T_q W^s(p,f)=l$.

For any $v \in T_q W^s(p,f)$, the norm of $Df^n v$
 with respect to the fixed Riemannian metric on $M$
 converges to zero as $n$ goes to infinity.
This implies that $\|D(Q^h \circ f^n) v\|$ converges to zero
 for $v \in T_q W^s(p,f)$.
For $w \in C(q,\theta,Q^h,P^h)$,
 we have $\|D(Q^h \circ f^n) w\| \geq \|DQ^h w\|$ for any $n \geq 0$
 by the cone condition.
Therefore, $T_q W^s(p,f) \cap C(q,\theta,Q^h,P^h)=\{0\}$.
Since $T_q W^s(p,f)$ is $l$-dimensional
 and the cone $C(q,\theta,Q^h,P^h)$ contains
 a $m$-dimensional subspace $\Ker DP^h$,
 we have a splitting $T_q M=T_q W^s(p,f) \oplus \Ker DP^h$.
In particular,
 any vector in $T_q M$ can be written as a sum of vectors
 in $T_q W^s(p,f)$ and $C(q,\theta,Q^h,P^h)$.

Fix $v_* \in \xi$ and we show $v_* \in T_q W^s(p,f)$.
Since $\xi$ is an element of $\Lambda^s(\wh{K},\wh{f})$,
 $Df^n(\xi)$ is a subset of $C(f^n(q),\theta,P^h,Q^h)$ for any $n \geq 1$.
This implies that $Df^n v_*$ is an element of $C(f^n(q),\theta,P^h,Q^h)$,
 and hence, $\|D(Q^h \circ f^n)v_*\| \leq \theta\|D(P^h \circ f^n) v_*\|$.
By the $(\theta,1,1)$-cone condition on $K$,
 there exists $\eta>1$ such that 
$\|DP^h v_*\| \geq \eta^n\|D(P^h \circ f^n) v_*\|$
 for any $n \geq 0$.
In particular,
 $\|D(Q^h \circ f^n) v_*\|$ converges to zero as $n$ goes to infinity.
Take $v \in T_q W^s(p,f)$ and $w \in C(q,\theta,Q^h,P^h)$
 such that $v_*=v+w$.
Then,
\begin{align*}
 \|DQ^h w\| & \leq \|D(Q^h \circ f^n) w\|
 \leq \|D(Q^h \circ f^n)v_*\|+\|D(Q^h \circ f^n) v\|.
\end{align*}
The right-hand side converges to zero as $n$ goes to infinity,
 and hence, $\|DQ^h w\|=0$.
Since $w$ is contained in $C(q,\theta,Q^h,P^h)$, this implies that $w=0$.
Therefore, $v_*$ is an element of $T_q W^s(p,f)$.
\end{proof}

\subsection{Blending machines for the lift to the Grassmannian bundle}
\label{sec:horseshoe}
We construct a diffeomorphism whose lift to the Grassmannian bundle
 admits a blending machine with respect to a natural splitting
 by following Barrientos and Raibekas' construction in \cite{BR17}.

For an $(l,m)$-matrix $A \in \Mat(l,m)$,
 let $\|A\|$ be the operator norm of $A \in \Mat(l,m)$
 with respect to the box norms on $\RR^l$ and $\RR^m$
 as a linear map from $\RR^m$ to $\RR^l$.
For $A \in \Mat(l,m)$ and $r>0$, we put
\begin{equation*}
 B^{(l,m)}(A,r)=\{A' \in \Mat(l,m) \mid \|A'-A\| \leq r\}.
\end{equation*}
 and write $B^{(l,m)}(r)$ when $A$ is the zero matrix.

Fix a quadruple $\vm=(m_1,m_2,m_3,m_4)$ of positive integers
 and put $|\vm|=m_1+m_2+m_3+m_4$.
Define an $(m_1+m_2,m_3+m_4)$-splitting $(P^h,Q^h) $ of $\RR^{|\vm|}$ by
\begin{alignat*}{4}
 P^h(x_1,x_2,x_3,x_4) & =(x_1,x_2), & \quad
 Q^h(x_1,x_2,x_3,x_4) & =(x_3,x_4)
\end{alignat*}
 for $(x_1,x_2,x_3,x_4) \in \RR^{|\vm|}$ with $x_i \in \RR^{m_i}$.
Remark that $\Ker DP^h=\{0\} \oplus \RR^{m_3+m_4}$ and
 $\Ker DQ^h=\RR^{m_1+m_2} \oplus \{0\}$.
Let $\wh{M}$ be the subset of $\Gr(\RR^{|\vm|},m_1+m_2)$ given by
\begin{equation*}
 \wh{M}=\{\xi \in \Gr(\RR^{|\vm|},m_1+m_2) \mid
  \xi \text{ is transverse to } \Ker DP^h=\{0\} \oplus \RR^{m_3+m_4}\}.
\end{equation*}
We define a map $\Pi:\wh{M} \ra \Mat(m_1+m_2+,m_3+m_4)$ by
\begin{equation*}
 D(P^h,Q^h)(\xi)=\{(w,\Pi(\xi)w) \mid w \in \RR^{m_1+m_2}\}.
\end{equation*}
In other words, $\Pi(\xi)$ is the unique element of
 $\Mat(m_1+m_2,m_3+m_4)$
 such that $DQ^h v=\Pi(\xi)(DP^h v)$ for any $v \in \xi$.
The map $(\pi, \Pi):\xi \mapsto (\pi(\xi),\Pi(\xi))$ is a diffeomorphism
 from $\wh{M}$ to $\RR^{|\vm|} \times \Mat(m_1+m_2,m_3+m_4)$.
For $i=1,2,3,4$, let $\pi_i:\wh{M} \ra \RR^{m_i}$
 be the $i$-th component of the natural projection
 $\pi:\wh{M} \ra
 \RR^{|\vm|}=\RR^{m_1} \times \RR^{m_2} \times \RR^{m_3} \times \RR^{m_4}$.
We define $\Pi_{ij}:\wh{M} \ra \Mat(m_i,m_j)$
 for $i=3,4$ and $j=1,2$ by
\begin{equation*}
 \Pi(\xi)=
\begin{bmatrix}
 \Pi_{31}(\xi) & \Pi_{32}(\xi)\\
 \Pi_{41}(\xi) & \Pi_{42}(\xi)
\end{bmatrix}.
\end{equation*}
We define a subset $\wh{M}(1)$ of $\wh{M}$ by
\begin{align*}
 \wh{M}(1) & =\{\xi \in \wh{M} \mid \|\Pi(\xi)\| \leq 1\}\\
 & =\{\xi \in \Gr(\RR^{|\vm|},m_1+m_2) \mid
 \xi \subset C(\pi(\xi),1,P^h,Q^h)\}.
\end{align*}
\begin{rmk}
\label{rmk:Pi norm}
For $A_{ij} \in B^{(m_i,m_j)}(r_{ij})$ $(i=3,4,j=1,2)$,
\begin{equation*}
 \left\|
\begin{bmatrix} A_{31} & A_{32} \\ A_{41} & A_{42}\end{bmatrix}
\right\| \leq
 \|A_{31}\|+ \|A_{32}\|+ \|A_{41}\|+ \|A_{42}\|
 \leq r_{31}+r_{32}+r_{41}+r_{42}.
\end{equation*}
\end{rmk}
Put
 $\wh{l} = m_1+m_2+m_3+m_3m_2$,
 $\wh{m} = m_4++m_3m_1+m_4m_1+m_4m_2$ and 
 let $E^+$ and $E^+$
 be $\wh{l}$- and $\wh{m}$-dimensional vector spaces given by
\begin{align}
\label{eqn:E-}
E^- &
 = \RR^{m_1} \times \RR^{m_2} \times \RR^{m_3} \times \Mat(m_3,m_2),\\
\label{eqn:E+}
E^+ &
 =  \RR^{m_4} \times \Mat(m_3,m_1)\times \Mat(m_4,m_1)
 \times \Mat(m_4,m_2).
\end{align}
We define a $(\wh{l},\wh{m})$-splitting
 $(\mc{P},\mc{Q}):\wh{M} \ra E^- \times E^+$ by
\begin{alignat}{4}
\label{eqn:mcPQ}
 \mc{P} & = (\pi_1,\pi_2,\pi_3,\Pi_{32},),& \quad
 \mc{Q} & = (\pi_4,\Pi_{31},\Pi_{41},\Pi_{42}).
\end{alignat}

Define a subset $\mc{Z}^+$ of $E^+$ by
\begin{align}
\label{eqn:Z+}
\mc{Z}^+ & =B^{m_4}(1) \times B^{(m_3,m_1)}(1/8)
 \times B^{(m_4,m_1)}(1/8) \times B^{(m_4,m_2)}(1/8).
\end{align}
For $\mu>0$, we also define subsets $\mc{Z}^-(\mu)$ and $\mc{Z}(\mu)$
 of $E^-$ and $\wh{M}$ respectively by
\begin{align}
\label{eqn:Z-}
\mc{Z}^-(\mu) & =B^{m_1}(1) \times B^{m_2}(1) \times B^{m_3}(1)
 \times B^{(m_3,m_2)}(1/4\mu),\\
\label{eqn:Z}
\mc{Z}(\mu) & =\{\xi \in \wh{M} \mid \mc{P}(\xi) \in \mc{Z}^-,
 \mc{Q}(\xi) \in \mc{Z}^+\}.
\end{align}
In the rest of this subsection,
 we will construct a diffeomorphism of $\RR^{|\vm|}$
 with a hyperbolic invariant set
 such that the lift of the diffeomorphism to the Grassmannian bundle
 admits a blending machine with respect to
 the splitting $(\mc{P},\mc{Q})$.
\begin{prop}
\label{prop:BR} 
Fix a quadruple $\vm=(m_1,m_2,m_3,m_4)$ of positive integers,
 positive constants $\lambda_*,\mu_*,\mu$ 
 with $\sqrt{\lambda_*}>\mu_*>\mu>2$,
 and open subsets $U^-, U^+$ of $\Int B^{m_1}(1)$, $\Int B^{m_4}(1)$.
Then, there exist a $C^\infty$ diffeomorphism $f$ of $\RR^{|\vm|}$,
 a finite set $I$,
 families  $(K_i)_{i \in I}$ and $\mc{B}=(\mc{R}_i)_{i \in I}$
 of mutually disjoint compact subsets of $\RR^{|\vm|}$
 and $\Gr(\RR^{|\vm|},m_1+m_2)$
 respectively which satisfy the following properties:
\begin{enumerate}
\item For any $i \in I$,
\begin{alignat*}{4}
P^h(K_i) & = B^{m_1}(4) \times B^{m_2}(4), &\hsp
Q^h(K_i) & \subset \Int B^{m_3}(4)  \times U^+,\\
Q^h(f(K_i)) & = B^{m_3}(4) \times B^{m_4}(4) & \hsp
P^h(f(K_i)) & \subset U^- \times \Int B^{m_2}(4).
\end{alignat*}
\item The restriction of the map $f$ to each $K_i$
 satisfies the $(1,1,1)$-cone condition
 with respect to $(P^h,Q^h)$.
\item The lift $\wh{f}$ of $f$ to the Grassmannian bundle
 $\Gr(\RR^{|\vm|},m_1+m_2)$
 satisfies the $(\theta,\mu_*^{-2},\lambda_*)$-cone condition on
 $\bigcup_{i \in I}\mc{R}_i$ for any $\theta>0$.
\item $\pi(\mc{R}_i) \subset \Int K_i$ and
 $\mc{R}_i \cup \wh{f}(\mc{R}_i) \subset \Int \wh{M}(1)$ for each $i \in I$. 
\item The family $\mc{B}$ is a blending machine for $\wh{f}$
 with respect to the splitting $(\mc{P},\mc{Q})$
 whose superposition region is $\mc{Z}(\mu)$.
\end{enumerate}
\end{prop}
As mentioned above, the proof of the proposition
 follows the construction by Barrientos and Raibekas in \cite{BR17}
 with minor adaptation to our setting.
Fix a quadruple $\vm=(m_1,m_2,m_3,m_4)$,
 constants $\lambda_*$, $\mu_*$, $\mu$ with $\sqrt{\lambda_*}>\mu_*>\mu>2$
 and open sets $U^-$, $U^+$ as in the proposition.
There exist a finite set $I$,
 and families $(q_i)_{i \in I}$, $(C_i)_{i \in I}$
 of points in $B^{m_3}(1)$, and $B^{m_3,m_2}(1/4\mu)$ respectively
 such that the families $(\Int B^{m_3}(q_i,1/4\mu))_{i \in I}$
 and  $(\Int B^{(m_3,m_2)}(C_i,1/8\mu^3))$ are
 open covers of $B^{m_3}(1)$ and $B^{(m_3,m_3)}(1/4\mu)$ respectively.
We take $\lambda>\lambda_*$ and
 families $(p_i^-)_{i \in I}$, $(p^+_i)_{i \in I}$
 of points of $\Int B^{m_1}(1)$, $\Int B^{m_4}(1)$ such that
 $(B^{m_1}(p^-_i,4/\lambda))_{i \in I}$
 and $(B^{m_4}(p^+_i,4/\lambda))_{i \in I}$ are families
 of mutually disjoint subsets of $U^-$ and $U^+$ respectively.
Put
\begin{align*}
 K_i^c & = \{(x_2,x_2) \in \RR^{m_2} \times \RR^{m_3} \mid
 x_2 \in B^{m_2}(4), x_3 \in B^{m_3}(q_i+C_ix_2,4/\mu)\},\\
 K_i & = B^{m_1}(4) \times K_i^c \times B^{m_4}(p^+_i,4/\lambda).
\end{align*}
 and define an affine diffeomorphism $f_i$ of $\RR^{|\vm|}$ by
\begin{align*}
 f_i(x_1,x_2,x_3,x_4)
 & =(\lambda^{-1}x_1+p^-_i, \mu^{-1}x_2,
  \mu(x_3-(q_i+C_ix_2)),\lambda(x_4-p^+_i)).
\end{align*}
Then, we have
\begin{align*}
 K_i & \subset B^{m_1}(4) \times B^{m_2}(4) \times B^{m_3}(1+(5/\mu))
 \times B^{m_4}(p^+_i,4/\lambda),\\
 f_i(K_i) & = B^{m_1}(p^-_i,4\lambda^{-1}) \times B^{m_2}(4\mu^{-1})
 \times B^{m_3}(4) \times B^{m_4}(4).
\end{align*}
In particular,
\begin{alignat*}{4}
 P^h(K_i) & = B^{m_1}(4) \times B^{m_2}(4), & \hsp
 Q^h(K_i) & \subset \Int B^{m_3}(4) \times U^+, \\
 Q^h(f_i(K_i)) & = B^{m_3}(4) \times B^{m_4}(4), & \hsp
 P^h(f_i(K_i)) & \subset U^- \times \Int B^{m_2}(4).
\end{alignat*}
Hence, $K_i$ satisfies the first property in the proposition.
\begin{figure}
\begin{center}
 \includegraphics[scale=0.6]{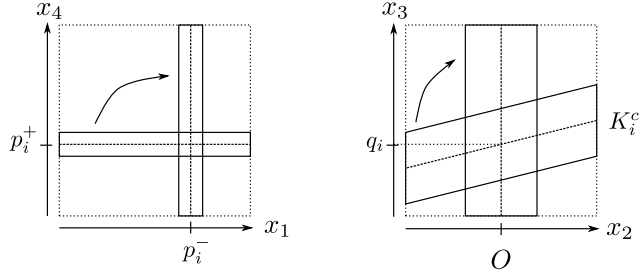} 
\end{center}
 \caption{the map $f$ on $K_i=B^{m_1}(4) \times K_i^c \times B^{m_4}(p^+_i4/\mu)$}
\label{fig:BR}
\end{figure}
Since both $(B^{m_1}(p^-_i,4/\lambda))_{i \in I}$
 and $(B^{m_4}(p^+_i,4/\lambda))_{i \in I}$
 are families of mutually disjoint sets in $\RR^{m_1}$ and $\RR^{m_4}$,
 both $(K_i)_{i \in I}$ and $(f_i(K_i))_{i \in I}$ are also families
 of mutually disjoint sets which are diffeomorphic to
 the unit square $[0,1]^{2m}$.
By lemma \ref{lemma:paste},
 there exists a $C^\infty$ diffeomorphism $f$ of $\RR^{|\vm|}$ such
 that $f=f_i$ on a neighborhood of each $K_i$.
Remark that $f$ satisfies the first condition in Proposition \ref{prop:BR}.
The next lemma shows that $f$ satisfies the second one.
\begin{lemma}
\label{lemma:BR hyp}
The map $f$ satisfies the $(1,1,1)$-cone condition
 with respect to $(P^h,Q^h)$ on each $K_i$.
\end{lemma}
\begin{proof}
Take $p \in K_i$ and
 $v=(v_1, v_2,v_3, v_4) \in T_p \RR^{|\vm|}$ with $v_i \in \RR^{m_i}$.
Then,
\begin{align*}
 \|D(P^h \circ f) v\| & = \|(\lambda^{-1}v_1,\mu^{-1}v_2)\|
 \leq \mu^{-1}\|(v_1,v_2)\|=\mu^{-1}\|DP^h v\|, \\
 \|D(Q^h \circ f) v\| & =\|(-C_iv_2+\mu v_3,\lambda v_4)\|
  \geq \mu\|(v_3,v_4)\|-\|C_i\|\|v_2\| \\
  & \geq \mu\|DQ^h v\|-(1/4\mu)\|DP^h v\| \geq \mu\|DQ^h v\|-\|DP^h v\|.
\end{align*}
If $v$ is contained in $C(p,1,Q^h,P^h)$ {\it i.e.}, $\|DQ^h v\| \geq \|DP^h v\|$,
 then 
\begin{align*}
 \|D(Q^h \circ f) v\| & \geq (\mu-1)\|DQ^h v\|
 \geq  (\mu-1)\|DP^h v\| \geq \mu(\mu-1)\|D(P^h \circ f)v\|.
\end{align*}
Hence,
 $\|D(Q^h \circ f) v\| \geq (\mu-1) \|DQ^h v\|$
 and $Df v$ is contained in $C(f(p),(\mu(\mu-1))^{-1},Q^h,P^h)$.
Remark that $\mu-1>1$ and $0<(\mu(\mu-1))^{-1}<1$ since $\mu>2$.
Similarly, for
 $w=(w_1,w_2,w_3,w_4) \in T_{f(p)} \RR^{|\vm|}$,
 we have
\begin{align*}
 \|D(P^h \circ f^{-1}) w\| & = \|(\lambda w_1,\mu w_2)\|
 \geq \mu\|(w_1,w_2)\|=\mu\|DP^h w\|, \\
 \|D(Q^h \circ f^{-1}) w\| & =\|(C_i w_2+\mu^{-1} w_3,\lambda^{-1} w_4)\|
  \leq \mu^{-1}\|(w_3,w_4)\|+\|C_i\|\|v_2\| \\
 & \leq \mu^{-1}\|DQ^h w\|+(4\mu)^{-1}\|DP^h w\|
  \leq \mu^{-1}(\|DQ^h w\|+\|DP^h w\|).
\end{align*}
If $w$ is contained in $C(f(p),1,P^h,Q^h)$, then 
\begin{align*}
 \|D(Q^h \circ f^{-1}) w\| &
 \leq 2\mu^{-1}\|DP^h w\| \leq 2\mu^{-2}\|D(P^h \circ f^{-1}) w\|.
\end{align*}
Hence, $\|D(P^h \circ f^{-1} w)\| \geq \mu\|DP^h w\|$ and
 $Df^{-1} w$ is contained in $C(f(p),2\mu^{-2},P^h,Q^h)$.
Remark that $\mu>2$, and hence, $0<2\mu^{-2}<1$.
Therefore, $f$ satisfies the $(1,1,1)$-cone condition for $(P^h,Q^h)$.
\end{proof}
By Proposition \ref{prop:transitive},
 $\Lambda(\bigcup_{i \in I}K_i,f)$
 is a locally maximal and topologically transitive hyperbolic set of $f$
 whose unstable index is $m_3+m_3$.
\begin{rmk}
\label{rmk:blender K}
Let $(P,Q)$ be the natural $(m_1+m_2+m_3,m_4)$-splitting
 of $\RR^{|\vm|}=\RR^{m_1+m_2+m_3} \times \RR^{m_4}$
 and define a family $(R_i)_{i \in I}$ of compact subsets of $\RR^{|\vm|}$ by
\begin{equation*}
 R_i=B^{m_1}(2) \times R^c_i \times B^{m_4}(p^+_i,\lambda^{-1}) 
 =\pi(\mc{R}_i)
\end{equation*}
 where $R^c_i$ and $\mc{R}_i$ are compact subsets of $\RR^{m_2+m_3}$
 and $\wh{M}$ defined by (\ref{eqn:Rc}) and (\ref{eqn:R}) below.
Then, one can check that the family $(R_i)_{i \in I}$
 is a blending machine for $f$ with respect to $(P,Q)$
 such that the superposition domain is $B^{|\vm|}(1)=[-1,1]^{|\vm|}$.
This implies that the hyperbolic set $\Lambda(\bigcup_{i \in I}K_i,f)$
 is a $cu$-blender of $uu$-index $m_4$
 by Proposition \ref{prop:blender}.
However, we do not need this fact to prove Theorem B.
\end{rmk}
Next, we construct a blending machine for the lift of $f$ with respect to
 the splitting  $(\mc{P},\mc{Q})$ of the subspace $\wh{M}$
 of the Grassmannian bundle $\Gr(\RR^{|\vm|},m_1+m_2)$.
Let $\wh{f}:\Gr(\RR^{|\vm|},m_-) \ra \Gr(\RR^{|\vm|},m_-)$
 be the lift of $f$.
Since $f$ preserves the vertical foliation
 $(\{x\} \times \RR^{m_3+m_4})_{x \in \RR^{m_1+m_2}}$,
 the restriction of the lift $\wh{f}$ to $\wh{M}$
 is a diffeomorphism of $\wh{M}$.
Put
\begin{equation*}
\wh{K}_i= \pi^{-1}(K_i) \cap \wh{M}(1)
 =\{\xi \in \wh{M} \mid \pi(\xi) \in K_i, \|\Pi(\xi)\| \leq 1\}
\end{equation*}
 and $\wh{K}=\bigcup_{i \in I}\wh{K}_i$.
 and define affine maps $F^-_i:E^- \ra E^-$ and $F^+_i:E^+ \ra E^+$ by 
\begin{align*}
F^-_i(x_1,x_2,x_3,A_{32})
 & =(\lambda^{-1} x_1+p^-_i, \mu^{-1} x_2,
 \mu(x_3-q_i-C_ix_2), \mu^2(A_{31}-C_i)),\\
F^+_i(x_4, A_{31}, A_{41}, A_{42})
 & = (\lambda(x_4-p^+_i), \lambda\mu A_{31},
 \lambda^2 A_{41}, \lambda\mu A_{42}).
\end{align*}
Then, by direct computation, we have
\begin{align*}
(\mc{P},\mc{Q}) \circ \wh{f} \circ (\mc{P},\mc{Q})^{-1}(\xi^-,\xi^+)
 & =(F^-_i(\xi^-),F^+_i(\xi^+))
\end{align*}
 for $(\xi^-,\xi^+) \in (\mc{P},\mc{Q})(\wh{K}_i)$.
\begin{lemma}
\label{lemma:BR cone}
For any $\theta>0$,
 the map $\wh{f}$ satisfies the
 $(\theta,\mu_*^{-2},\lambda_*)$-cone condition with respect
 to $(\mc{P},\mc{Q})$ on $\wh{K}$.
\end{lemma}
\begin{proof}
For $\xi \in \wh{K}$ with $\pi(\xi) \in K_i$
 and $v \in T_\xi \wh{M}$, we have
\begin{alignat*}{4}
D(\mc{P} \circ \wh{f})(v) & =DF^-_i(D\mc{P}v) , &\hsp
D(\mc{Q} \circ \wh{f})(v) & =DF^+_i(D\mc{Q}v).
\end{alignat*}
It is easy to see that 
\begin{align*}
\|DF^-_i v \| & \leq \max\{\lambda^{-1},\mu^{-1},\mu(1+\|C_i\|),\mu^2\} \|v\|
 \leq \mu^2 \|v\|,\\
\|DF^+_i w\| & \geq  \min\{\lambda,\lambda\mu,\lambda^2\}\|w\|
 \geq \lambda\|w\|
\end{align*}
 for any $v \in T_{\mc{P}(\xi)}E^-$ and $w \in T_{\mc{Q}(\xi)} E^+$.
Since $\lambda>\lambda_*>\mu_*^2>\mu^2$,
 this implies that $\wh{f}$ satisfies the
 $(\theta,\mu_*^{-2},\lambda_*)$-cone condition with respect
 to $(\mc{P},\mc{Q})$ on $\wh{K}$ for any $\theta>0$.
\end{proof}

Next, we define a blending machine $\mc{B}=(\mc{R})_{i \in I}$ for $\wh{f}$
 with respect to $(\mc{P},\mc{Q})$ and show that it satisfies
 the desired properties.
For each $i \in I$, we define subsets $R_i^c$, $\mc{R}^-_i$, $\mc{R}^+_i$,
 and $\mc{R}$ of $\RR^{m_2+m_3}$, $E^-$, $E^+$, and $\wh{M}$ respectively by
\begin{align}
\label{eqn:Rc}
R^c_i & =\{(x_2,x_3) \in \RR^{m_2} \times \RR^{m_3} \mid
 x_1 \in B^{m_2}(2), B^{m_3}(q_i+C_i x_2,3/4\mu)\},\\
\label{eqn:R-}
\mc{R}^-_i & = B^{m_1}(2) \times R^c_i
 \times B^{(m_3,m_2)}(C_i,1/8\mu^{3}),\\
\label{eqn:R+}
\mc{R}^+_i & = B^{m_4}(p^+_i,1/\lambda) \times
 B^{(m_3,m_1)}(1/8\lambda\mu) \\
 & \quad  \times B^{(m_4,m_1)}(1/8\lambda^2) \times
 B^{(m_4,m_2)}(1/8\lambda\mu), \nonumber\\
\label{eqn:R}
\mc{R}_i & = (\mc{P},\mc{Q})^{-1}(\mc{R}_i^- \times \mc{R}_i^+)
 =\{\xi \in \wh{M} \mid
  \mc{P}(\xi) \in \mc{R}^-_i, \mc{Q}(\xi) \in \mc{R}^+_i\}.
\end{align}
By direct computation, we have
\begin{align}
\label{eqn:fR-}
 F^-_i(\mc{R}^-_i) 
 & =B^{m_1}(p^-_i,2/\lambda)
 \times B^{m_2}(2/\mu) \times B^{m_3}(3/4) \times B^{(m_3,m_2)}(1/8\mu),\\
\label{eqn:fR+}
F^+_i(\mc{R}^+_i)
 & = B^{m_4}(1) \times B^{(m_3,m_1)}(1/8) \times 
 B^{(m_4,m_1)}(1/8) \times  B^{(m_4,m_2)}(1/8).
\end{align}
It is easy to see
 that the compact sets $\mc{R}_i$ and $\wh{f}(\mc{R}_i)$ are
 $(\mc{P},\mc{Q})$-rectangles.
By Remark \ref{rmk:Pi norm},
 they are contained in
 $\Int \wh{M}(1)=\Int \{\xi \in \wh{M} \mid \|\Pi(\xi)\| \leq 1\}$.
Since $(\mc{P},\mc{Q}) \circ \wh{f}=(\mc{P},F^+_i \circ \mc{Q})$
 is a diffeomorphism between $\mc{R}_i$
 and $\mc{R}_i^- \times F_i^+(\mc{R}_i^+)$,
 the compact set $\mc{R}_i$ is a $(\mc{P},\mc{Q} \circ \wh{f})$-rectangle.

The following lemma finishes the proof of Proposition \ref{prop:BR}.
\begin{lemma}
The family $\mc{B}=(\mc{R}_i)_{i \in I}$ is a blending machine
 for $\wh{f}$ with respect to the splitting $(\mc{P},\mc{Q})$
 such that its superposition domain is $\mc{Z}(\mu)$.
\end{lemma}
\begin{proof}
By the equations (\ref{eqn:Z+}), (\ref{eqn:Z-}),
 (\ref{eqn:R-}), (\ref{eqn:fR-}), and (\ref{eqn:fR+}),
\begin{gather*}
\mc{Q}(\mc{R}_i)  = \mc{R}^+_i \subset \Int \mc{Z}^+
 = \Int \mc{Q}(\mc{Z}(\mu)),\\
\mc{P}(\wh{f}(\mc{R}_i))  = F^-_i(\mc{R}^-_i) \subset \Int \mc{Z}^-(\mu)
 = \Int \mc{P}(\mc{Z}(\mu)), \\
\mc{Q}(\wh{f}(\mc{R}_i))  = F^+_i(\mc{R}^+_i)=
 \mc{Z}^+=\mc{Q}(\mc{Z}(\mu)).
\end{gather*}
Therefore, the family $(\mc{R}_i)_{i \in I}$ satisfies the second condition
 in the definition of a blending machine.
Recall that the families $(B^{m_3}(q_i,1/4\mu))_{i \in I}$
 and $(B^{(m_3,m_2)}(C_i,1/8\mu^3))_{i \in I}$
 are open covers of $B^{m_3}(1)$ and $B^{(m_3,m_2)}(1/4\mu)$ respectively.
Since $\|C_i\| \leq 1/4\mu$,
 the set $R_i^c$ contains $B^{m_2}(1) \times B^{m_3}(q_i,1/4\mu)$ for each $i$.
This implies that the family $(\Int\mc{P}(\mc{R}_i))_{i \in I}$
 of open subsets of $E^-$ covers $\mc{P}(\mc{Z}(\mu))=\mc{Z}^-(\mu)$.
By Lemma \ref{lemma:BR cone},
 the map $\wh{f}$ satisfies the $(\theta,\mu_*^{-2},\lambda_*)$-cone
 condition for any $\theta>0$.
For small $\theta>0$,
 the Lebesgue number of the open cover
 $(\Int \mc{P}(\mc{R}_i))_{i \in I}$ of $\mc{P}(\mc{Z}(\mu))$
 is greater than $\theta \diam \mc{Q}(\mc{Z})(\mu)$.
Therefore, the family $(\mc{R}_i)_{i \in I}$ is a blending machine
 for $\wh{f}$ with respect to $(\mc{P},\mc{Q})$
 whose superposition domain is $\mc{Z}(\mu)$.
\end{proof}
\begin{rmk}
\label{rmk:blender Gr}
Define a $(m_1+m_2,(m_1+m_2+1)(m_3+m_4))$-splitting $(\mc{P}^h,\mc{Q}^h)$
 of $\wh{M}$ by $\mc{P}^h=(\pi_1,\pi_2)$
 and $\mc{Q}^h=(\pi_3,\pi_4,\Pi)$.
Then, it is possible to show that
 $\wh{f}$ satisfies the $(1,1,1)$-cone condition on $\wh{M}$
 with respect to $(\mc{P}^h,\wh{Q}^h)$
 and $\Lambda(\wh{K},\wh{f})$ is a locally maximal hyperbolic set
 of the unstable index $(m_1+m_2+1)(m_3+m_4)$.
Since $\bigcup_{i \in I}\mc{R}_i$ is contained in $\Int \wh{K}$
 and $\mc{B}=(\mc{R}_i)_{i \in I}$ is a blender machine,
 the hyperbolic set $\Lambda(\wh{K},\wh{f})$ is
 a $cu$-blender of $uu$-index $\wh{m}$.
However, we do not need this fact to prove Theorem B.
\end{rmk}

\subsection{A map connecting blending machines}
\label{sec:connecting}
Put
\begin{align*}
 \vm_1 & =(m_{1,1},m_{1,2},m_{1,3},m_{1,4})=(1,m-1,m-1,1),\\
 \vm_2 & =(m_{2,1},m_{2,2},m_{2,3},m_{2,4})=(1,m-1,1,m-1).
\end{align*}
In the next subsection,
 we will prove Theorem B by applying Theorem \ref{thm:intersection}
 to two blending machines in the Grassmannian bundle
 obtained from Proposition \ref{prop:BR} with $\vm=\vm_1, \vm_2$.
In this subsection, we construct a map $f_\sh$ 
 whose lift $\wh{f}_\sh$ provides
 intersection of the forward invariant sets of the blending machines.

As before, let $(P^h,Q^h)$ be the natural $(m,m)$-splitting
 of $\RR^{2m}=\RR^m \times \RR^m$
 and $\wh{M}$ the open subset of $\Gr(\RR^{2m},m)$
 which consists of $m$-dimensional planes transverse to $\{0\} \oplus \RR^m$.
For $\tau=1,2$, let $E^-_\tau$ and $E^+_\tau$ be
 the linear spaces defined by (\ref{eqn:E-}), (\ref{eqn:E+}),
 and $(\mc{P}_\tau,\mc{Q}_\tau):\wh{M} \ra E^-_\tau \oplus E^+_\tau$
 the splitting of $\wh{M}$ defined by (\ref{eqn:mcPQ})
 associated with the quadruple $\vm_\tau$.
Remark that $(\mc{P}_1,\mc{Q}_1)$ is a $(m^2,2m)$-splitting
 and $(\mc{P}_2,\mc{Q}_2)$ is a $(2m,m^2)$-splitting
 of the $(m^2+2m)$-dimensional manifold $\wh{M}$.
\begin{prop}
\label{prop:h-sh}
For any positive numbers $\theta,\nu$,
 any open neighborhoods $\mc{V}^-_1$, $\mc{V}^-_2$
 of the origins of $E^-_1$, $E^-_2$,
 and any non-empty compact subsets $\mc{K}^+_1, \mc{K}^+_2$
 of $E^+_1, E^+_2$ respectively,
 there exists a pair $(f_\sh,\mc{R}_\sh)$
 of a $C^\infty$ diffeomorphism of $\RR^{2m}$
 and a compact subset of $\wh{M}$ which satisfies
 the following properties:
\begin{enumerate}
 \item $f_\sh(B^m(1/4) \times B^m(4)) \subset B^m(1/4) \times \RR^m$.
 \item The map $f_\sh$ lifts to a diffeomorphism of $\wh{M}$
 which satisfies the $(\theta,\nu,\nu)$-cone condition
 with respect to the splittings $(\mc{Q}_2,\mc{P}_2)$
 and $(\mc{P}_1,\mc{Q}_1)$.
 \item $\mc{R}_\sh$ is a $(\mc{Q}_2,\mc{Q}_1 \circ \wh{f}_\sh)$-rectangle
 such that
\begin{alignat*}{4}
 \mc{P}_2(\mc{R}_\sh) & \subset \mc{V}^-_2,
&\hsp
 \mc{Q}_2(\mc{R}_\sh) & = \mc{K}^+_2,
& \hsp
 \mc{P}_1(\wh{f}_\sh(\mc{R})_\sh)
  & \subset \mc{V}^-_1,
&\hsp
 \mc{Q}_1(\wh{f}_\sh(\mc{R}_\sh)) & = \mc{K}^+_1.
\end{alignat*}
\end{enumerate}
\end{prop}
We recommend the reader to compare the above conditions
 with the assumptions on $(R_\sh,f_\sh)$ in Theorem \ref{thm:intersection}.

Let $O$ be the origin of $\RR^{2m}$
 and $\wh{O}$ the $m$-dimensional linear subspace $\RR^m \times \{0\}$
 in $T_O\RR^{2m}$.
The linear subspace $\wh{O}$ belongs to $\wh{M}$
 as an element of $\Gr(\RR^{2m},m)$,
It satisfies that $\pi(\wh{O})=O$,
 $\mc{P}_\tau(\wh{O})=0$, and $\mc{Q}_\tau(\wh{O})=0$ for $\tau=1,2$.
The first step to construct the required map $f_\sh$
 is to find a diffeomorphism $h$
 such that the lift $\wh{h}$ is a diffeomorphism of $\wh{M}$
 with $\wh{h}(\wh{O})=\wh{O}$
 and the map $(\mc{Q}_2,\mc{Q}_1 \circ \wh{h})$
 is an embedding on a small neighborhood of $\wh{O}$.
The condition $\wh{h}(\wh{O})=\wh{O}$ means
 that $Dh$ preserves the subspace $\RR^m \times \{0\}$ of $T_O\RR^{2m}$.
In particular,
 $h(\RR^m \times \{0\})$ is tangent to $\RR^m \times {0}$
 at the origin.
If $(\mc{Q}_2,\mc{Q}_1 \circ \wh{h})$
 is an embedding on a small neighborhood of $\wh{O}$,
 then $D\wh{h}(\Ker(D\mc{Q}_2))$ is transverse to $\Ker \mc{Q}_1$
 in $T_O\RR^{2m}$.
This condition requires that the tangency at the origin is quadratic.
From these observations,
 we define a diffeomorphism $h$ of $\RR^{2m}$ by
\begin{equation*}
 h(y_1,\vec{y}_2,y_3,\vec{y}_4)
 =\left(y_1, \vec{y}_2, \vec{y}_4+y_1\vec{y}_2, y_3+\frac{y_1^2}{2} \right)
\end{equation*}
 for $(y_1,\vec{y}_2,y_3,\vec{y}_4) \in \RR^{2m}$
 with $y_1,y_3 \in \RR$ and $\vec{y}_2,\vec{y}_4 \in \RR^{m-1}$
 and show this map has the required properties above.
The map $h$ preserves the vertical plane
 $\{x\} \times \RR^m$ for each $x \in \RR^m$.
This implies that $h$ lifts to a diffeomorphism $\wh{h}$ of $\wh{M}$.
By direct computation, we have
\begin{align}
\label{eqn:wh-h}
\lefteqn{(\mc{P}_1,\mc{Q}_1) \circ \wh{h} \circ (\mc{P}_2,\mc{Q}_2)^{-1}
(y_1,\vec{y}_2,y_3,Y_{32},\vec{y}_4,y_{31},\vec{y}_{41},Y_{42})}\\
& =
\left(y_1,\vec{y}_2,\vec{y}_4+y_1\vec{y_2},Y_{42}+y_1 I,
 y_3+\frac{y_1^2}{2}, \vec{y}_{41}+\vec{y}_2,y_{31}+y_1,Y_{32}\right),
 \nonumber
\end{align}
 where $y_1, y_3,y_{31} \in \RR=\Mat(1,1)$,
 $\vec{y}_2, \vec{y}_4,\vec{y}_{41} \in \RR^{m-1}=\Mat(m-1,1)$,
 $Y_{32} \in \Mat(1,m-1)$, $Y_{42} \in \Mat(m-1,m-1)$,
 and $I$ is the unit matrix of size $(m-1)$.
\begin{lemma}
\label{lemma:h-sh}
The map $\wh{h}$ fixes the point $\wh{O}$.
There exist
 a compact neighborhood $\mc{V}_0$ of $\wh{O}$ in $\wh{M}$ 
 and constants $\alpha,\beta>0$ such that
 the restriction of $(\mc{Q}_2, \mc{Q}_1 \circ \wh{h})$ to $\mc{V}_0$
 is an embedding onto a compact neighborhood of the origin $(0,0)$
 in $E^+_2 \times E^-_1$ and
\begin{gather}
\label{eqn:h-sh cone 1}
\wh{h}(C(\xi,\alpha,\mc{P}_2,\mc{Q}_2))
 \subset C(\wh{h}(\xi),\beta,\mc{Q}_1,\mc{P}_1),\\
\label{eqn:h-sh cone 2}
 \wh{h}^{-1}(C(\wh{h}(\xi),\alpha,\mc{P}_1,\mc{Q}_1))
 \subset C(\xi,\beta,\mc{Q}_2,\mc{Q}_2)
\end{gather}
 for any $\xi \in \mc{V}_0$.
\end{lemma} 
\begin{proof}
The map $h$ fixes the origin $O$ of $\RR^{2m}$
 and the differential $Dh$ preserves
 the linear subspace $\wh{O}=\RR^m \times \{0\}$ of $T_O\RR^{2m}$.
This implies that $\wh{h}(\wh{O})=\wh{O}$.
By (\ref{eqn:wh-h}),
\begin{align*}
\mc{Q}_1 \circ \wh{h}_0 \circ (\mc{P}_2,\mc{Q}_2)^{-1}
(y_1,\vec{y}_2,y_3,Y_{32},0,0,0,0)
 =\left(y_3+\frac{y_1^2}{2},\vec{y}_2,y_1,Y_{32}\right)
\end{align*}
 for $(y_1,\vec{y}_2,y_3,Y_{32}) \in E^-_2$.
This implies that 
 the restriction of $D(\mc{Q}_1 \circ \wh{h}):T_{\wh{O}}\wh{M} \ra E^+_1$
 to $\Ker D\mc{Q}_2$ is an isomorphism.
In particular,
 $\Ker D\mc{Q}_2 \cap \Ker D(\mc{Q}_1 \circ \wh{h})=\{0\}$,
 and hence, $D(\mc{Q}_2,\mc{Q}_1 \circ \wh{h})$ is injective.
Since $\dim \wh{M}=\dim (E^+_2 \oplus E^+_1)=m^2+2m$,
 the differential $D(\mc{Q}_2,\mc{Q}_1 \circ \wh{h})$ is an isomorphism
 from $T_{\wh{O}}\wh{M}$ to $E^+_2 \oplus E^+_1$.
By the inverse function theorem,
 there exists an open neighborhood $\mc{V}'$ such that
 the restriction of $(\mc{Q}_2,\mc{Q}_1 \circ \wh{h})$ to $\mc{V}'$
 is an embedding onto an open neighborhood of the origin $(0,0)$
 of $E^+_2 \oplus E^+_1$.
Since $\Ker D\mc{Q}_2 \cap \Ker D(\mc{Q}_1 \circ \wh{h})=\{0\}$,
 we have
\begin{alignat*}{4}
  \Ker D\mc{Q}_2 \cap D\wh{h}^{-1}(\Ker D\mc{Q}_1) & = \{0\}, & \hsp
  D\wh{h}(\Ker D\mc{Q}_2) \cap \Ker D \mc{Q}_1 & =\{0\}.
\end{alignat*}
This implies that there exist $\alpha,\beta>0$ such that
\begin{align*}
 D\wh{h}^{-1}(C(\wh{O},\alpha,\mc{P}_1,\mc{Q}_1))
 & \subset C(\wh{O},\beta/2,\mc{Q}_2,\mc{P}_2), \\
 D\wh{h}(C(\wh{O},\alpha,\mc{P}_2,\mc{Q}_2))
 & \subset C(\wh{O},\beta/2,\mc{Q}_1,\mc{P}_1).
\end{align*}
We take a compact neighborhood $\mc{V}_0\subset \mc{V}'$ of $\wh{O}$
 such that
\begin{align*}
 D\wh{h}^{-1}(C(\wh{h}(\xi),\alpha,\mc{P}_1,\mc{Q}_1))
 & \subset C(\xi,\beta,\mc{Q}_2,\mc{P}_2), \\
 D\wh{h}(C(\xi,\alpha,\mc{P}_2,\mc{Q}_2))
 & \subset C(\wh{h}(\xi),\beta,\mc{Q}_1,\mc{P}_1).
\end{align*}
 for any $\xi \in \mc{V}_0$.
\end{proof}

We will construct the map $f_\sh$ in Proposition \ref{prop:h-sh}
 by composing $h$ with hyperbolic linear maps.
For $\tau=1,2$, we define linear isomorphisms
 $g_\tau:\RR^{2m} \ra \RR^{2m}$, $\mc{L}^-_\tau:E^-_\tau \ra E^-_\tau$,
 and $\mc{L}^+_\tau:E^+_\tau \ra E^+_\tau$ by
\begin{align*}
 g_\tau(z_1,z_2,z_3,z_4)
 & =(e^{-3}z_1,e^{-1}z_2,e^{-2}z_3,e z_4),\\
 \mc{L}^-_\tau(z_1,z_2,z_3,Z_{32})
 & = (e^{-3}z_1,e^{-1}z_2,e^{-2}z_3,e^{-1}Z_{32}),\\
 \mc{L}^+_\tau(z_4,Z_{31},Z_{41},Z_{42})
 & =(e z_4,e Z_{31},e^4Z_{41},e^2Z_{42})
\end{align*}
 where $z_i \in \RR^{m_{\tau,i}}$ and
 $Z_{jk} \in \Mat(m_{\tau,j},m_{\tau,k})$.
It is easy to check that
 $g_\tau$ lifts to a diffeomorphism $\wh{g}$ of $\wh{M}$
 and it satisfies that
\begin{equation*}
 (\mc{P}_\tau,\mc{Q}_\tau) \circ \wh{g}_\tau
 \circ (\mc{P}_\tau,\mc{Q}_\tau)^{-1}
 (\xi^-,\xi^+)=(\mc{L}^-_\tau(\xi^-),\mc{L}^+_\tau(\xi^+))
\end{equation*}
 for any $(\xi^-,\xi^+) \in E^-_\tau \times E^+_\tau$.
The linear maps $\mc{L}^-_\tau$ and $\mc{L}^+_\tau$ satisfy that
\begin{equation*}
 \|D\mc{L}^-_\tau\| = \|(D\mc{L}^+_\tau)^{-1}\|= e^{-1}<1.
\end{equation*}
This implies that
\begin{align*}
 \wh{g}_\tau^n(C(\xi,\eta,\mc{Q}_\tau,\mc{P}_\tau))
 &  \subset C(\wh{g}_\tau^n(\xi),e^{-2n}\eta,\mc{Q}_\tau,\mc{P}_\tau),\\
 \wh{g}_\tau^{-n}(C(\wh{g}_\tau^n(\xi),\eta,\mc{P}_\tau,\mc{Q}_\tau))
 & \subset C(\xi,e^{-2n}\eta,\mc{P}_\tau,\mc{Q}_\tau)
\end{align*}
 for any $n \geq 1$, $\xi \in \wh{M}$, and $\eta>0$.

For $n \geq 1$, we define a diffeomorphism $h_n$ of $\RR^{2m}$ by
\begin{align*}
h_n(y_1,\vec{y}_2,y_3,\vec{y}_4)
 & = g_1^n \circ h_0 \circ g_2^{-n}(y_1,\vec{y}_2,y_3,\vec{y}_4)\\
 & =\left(y_1,\vec{y}_2,e^{-3n}\vec{y}_4+e^{2n} y_1\vec{y}_2,
 e^{3n}y_3+\frac{e^{7n}}{2}y_1^2\right).
\end{align*}
 and put $\mc{R}_n=\wh{g}_2^{n}(\mc{R}'_n)$.
We will prove Proposition \ref{prop:h-sh}
 by showing that the pair $(h_n,\mc{R}_n)$
 satisfies the conditions required for $(f_\sh,\mc{R}_\sh)$
 in the proposition for some large $n$.
\begin{proof}
[Proof of Proposition \ref{prop:h-sh}]
Notice that $h_n(B^m(1/4) \times B^m(4)) \subset B^m(1/4) \times \RR^m$
 for any $n \geq 1$.
Fix constants $\theta,\nu>0$.
For $\tau=1,2$, we also fix
 an open neighborhood $\mc{V}^-_\tau$ of the origin in $E^-_\tau$
 and a non-empty compact subset $\mc{K}^+_\tau$ of $E^+_\tau$.
Let $\mc{V}_0$ and $\alpha,\beta$ be the compact neighborhood of $\wh{O}$
 and positive constants obtained in Lemma \ref{lemma:h-sh}.
Since $(\mc{L}^+_\tau)^{-1}$ is uniformly contracting 
 and $(\mc{Q}_2,\mc{Q}_1 \circ \wh{h})$ sends the neighborhood $\mc{V}_0$
 of $\wh{O}$ to a neighborhood of the origin $(0,0)$ of $E^+_2 \times E^+_1$,
 there exists $n_1 \geq 1$ such that
\begin{equation*}
 (\mc{L}^+_2)^{-n}(\mc{K}^+_2) \times
 (\mc{L}^+_1)^{-n}(\mc{K}^+_1) 
 \subset (\mc{Q}_2, \mc{Q}_1 \circ \wh{h})(\mc{V}_0).
\end{equation*}
 for any $n \geq n_1$.
Put
\begin{equation*}
 \mc{R}'_n=(\mc{Q}_2, \mc{Q}_1 \circ \wh{h})^{-1}
 ((\mc{L}^+_2)^{-n}(\mc{K}^+_2) \times
 (\mc{L}^+_1)^{-n}(\mc{K}^+_1))
\end{equation*}
 and $\mc{R}_n=\wh{g}_2^n(\mc{R}'_n)$.
Remark that $\mc{R}'_n$ is a $(\mc{Q}_2, \mc{Q}_1 \circ \wh{h})$-rectangle
 with $\mc{Q}_2(\mc{R}'_n)= (\mc{L}^+_2)^{-n}(\mc{K}^+_2)$
 and
 $\mc{Q}_1 \circ \wh{h}(\mc{R}'_n)=(\mc{L}^+_1)^{-n}(\mc{K}^+_1)$.
Since $\wh{h}_n=\wh{g}_1^n \circ \wh{h} \circ \wh{g}_2^{-n}$
 and $\mc{Q}_\tau \circ \wh{g}_\tau=\mc{L}^+_\tau \circ \mc{Q}_\tau$
 for $\tau=1,2$, we have
\begin{align*}
 (\mc{Q}_2,\mc{Q}_1 \circ \wh{h}_n)
 & =((\mc{L}^+_2)^n \circ \mc{Q}_2, (\mc{L}^+_1)^n \circ \mc{Q}_1 \circ \wh{h})
 \circ \wh{g}_2^{-n}\\
 & =((\mc{L}^+_2)^n \times (\mc{L}^+_1)^n) \circ (\mc{Q}_2,\mc{Q} \circ \wh{h})
 \circ \wh{g}_2^{-n}.
\end{align*}
This implies that $\mc{R}_n=\wh{g}_2^n(\mc{R}'_n)$
 is a $(\mc{Q}_2,\mc{Q}_1 \circ \wh{h}_n)$-rectangle with
\begin{alignat*}{4}
\mc{Q}_2(\mc{R}_n) & = (\mc{L}^+_2)^n(\mc{Q}_2(\mc{R}'_n))
 = \mc{K}^+_2, &\hsp
\mc{Q}_1(\wh{h}_n(\mc{R}_n))
 & = (\mc{L}^+_1)^n(\mc{Q}_1 \circ \wh{h}_0(\mc{R}'_n))
 = \mc{K}^+_1.
\end{alignat*}
Since $\mc{L}^-_\tau$ is uniformly contracting for each $\tau=1,2$,
 there exists $n_2 \geq n_1$ such that
\begin{alignat*}{4}
(\mc{L}^-_2)^n(\mc{P}_2(\mc{R}'_n))
 & \subset  \mc{V}^-_2, & \hsp
(\mc{L}^-_1)^n(\mc{P}_1 \circ \wh{h}(\mc{R}'_n))
 & \subset  \mc{V}^-_2
\end{alignat*}
 for any $n \geq n_2$.
Similar to the above, we have
\begin{alignat*}{4}
\mc{P}_2(\mc{R}_n) & = (\mc{L}^-_2)^n(\mc{P}_2(\mc{R}'_n))
 \subset \mc{V}^-_2, & \hsp
\mc{P}_1(\wh{h}_n(\mc{R}_n))
 & = (\mc{L}^-_1)^n(\mc{P}_1 \circ \wh{h}(\mc{R}'_n))
 \subset \mc{V}^-_1.
\end{alignat*}
By the conditions (\ref{eqn:h-sh cone 1}), (\ref{eqn:h-sh cone 2})
 on cones and the compactness of $\mc{V}_0$,
 there exists $\epsilon>0$ such that
\begin{alignat*}{4}
\|D(\mc{Q}_1 \circ \wh{h})v_2\| & \geq \epsilon \|D\mc{P}_2 v_2\|, & \hsp
\|D(\mc{Q}_2 \circ \wh{h}^{-1})v_1\| & \geq \epsilon \|D\mc{P}_1 v_1\|
\end{alignat*}
 for any $\xi \in \mc{V}_0$,
 $v_2 \in C(\xi,\alpha,\mc{P}_2,\mc{Q}_2)$,
 and $v_1 \in C(\wh{h}(\xi),\alpha,\mc{P}_1,\mc{Q}_1)$.
Take $n \geq n_2$ such that
$e^{-2n}\theta<\alpha$, $e^{-2n}\beta<\theta$,
 and $e^{2n}\epsilon>\nu$.
For any $\xi \in \mc{R}_n$,
 and $v \in C(\xi,\theta,\mc{P}_2,\mc{Q}_2)$, we have
\begin{alignat*}{4}
 D\wh{g}_2^{-n} v
 & \in  C(\wh{g}_2^{-n}(\xi),e^{-2n}\theta,\mc{P}_2,\mc{Q}_2),
 & \hsp
 \|D(\mc{P}_2 \circ \wh{g}_2^{-n}) v\|
 & \geq e^n\|D\mc{P}_2 v\|.
\end{alignat*}
Since $g_2^{-n}(\xi)$ is contained in $\mc{R}'_n \subset \mc{V}_0$
 and $e^{-2n}\theta<\alpha$,
 we also have
\begin{alignat*}{4}
 D(\wh{h} \circ \wh{g}_2^{-n}) v
 & \in C(\wh{h} \circ \wh{g_2}^{-n}(\xi),\beta,\mc{Q}_1,\mc{P}_1),
 & \hsp
 \|D(\mc{Q}_1 \circ \wh{h} \circ \wh{g}_2^{-n}) v\|
 & \geq \epsilon e^n\|D\mc{P}_2 v\|.
\end{alignat*}
In the same way,
\begin{alignat*}{2}
 D\wh{h}_n v
= D(\wh{g}_1^n \circ \wh{h}_0 \circ \wh{g}_2^{-n}) v
 & \in C(\wh{h}_n(\xi),e^{-2n}\beta,\mc{Q}_1,\mc{P}_1),
 & \hsp
 \|D\wh{h}_n v\|
 & \geq \epsilon e^{2n}\|D\mc{P}_2 v\|.
\end{alignat*}
Since $e^{-2n}\beta<\theta$ and $e^{2n}\epsilon>\nu$,
 the map $\wh{h}_n$ satisfies the $(\theta,\nu,\nu)$-cone condition
 with respect to the splittings $(\mc{Q}_2,\mc{P}_2)$
 and $(\mc{P}_1,\mc{Q}_1)$.
\end{proof}

\subsection{Robust tangency of the largest codimension}
\label{sec:largest}
Fix $m \geq 2$ and $\lambda_*,\mu_*,\mu$
 with $\sqrt{\lambda_*}>\mu_*>\mu>2$.
By $(P^h,Q^h)$,
 we denote the natural $(m,m)$-splitting of $\RR^{2m}$
 given by
 $P^h(x,y)= x$ and $Q^h(x,y)=y$ for $x,y \in \RR^m$.
Let $\wh{M}$ be the open subset of the Grassmannian bundle
 $\Gr(\RR^{2m},m)$
 which consists of $m$-dimensional linear subspaces
 transverse to $\Ker DP^h$.
As in the previous subsection, we put
\begin{alignat*}{4}
\vm_1 & = (1,m-1,m-1,1), & \hsp
 (\wh{l}_1,\wh{m}_1) & = (m^2,2m), \\
 \vm_2 & = (1,m-1,1,m-1), & \hsp
 (\wh{l}_2,\wh{m}_2) & = (2m,m^2).
\end{alignat*}
For each $\tau=1,2$, let $E^-_\tau$ and $E^+_\tau$ be the linear spaces
 defined by (\ref{eqn:E-}) and (\ref{eqn:E+}),
 $(\mc{P}_\tau,\mc{Q}_\tau):\wh{M} \ra E^-_\tau \times E^+_\tau$
 the $(\wh{l}_\tau,\wh{m}_\tau)$-splitting
 of $\wh{M}$ defined by (\ref{eqn:mcPQ}),
 and $\mc{Z}^-_\tau(\mu),\mc{Z}^+_\tau,\mc{Z_\tau}(\mu)$ be subsets
 of $E^-_\tau, E^+_\tau, \wh{M}$ defined by
 (\ref{eqn:Z-}),  (\ref{eqn:Z+}), (\ref{eqn:Z}) for $\vm=\vm_\tau$.
We put
\begin{alignat*}{6}
 U^-_1 & = U^+_1 = (-1,-1/2), &\hsp
 U^-_2 & = (1/2,1), &\hsp 
 U^+_2 & = (1/2,1)^{m-1}.
\end{alignat*}
For $\tau=1,2$,
 let $f_\tau$, $I_\tau$, $(K_{\tau,i})_{i \in I_\tau}$,
 $\mc{B}_\tau=(\mc{R}_{\tau,i})_{i \in I_\tau}$,
 be the diffeomorphism, families of subsets of
 $\RR^{2m}$ and $\wh{M}$ obtained in Proposition \ref{prop:BR}
 for the quadruple $\vm_\tau$,
 constants $\lambda_*,\mu_*,\mu$, and open sets $U^-_\tau$, $U^+_\tau$.
We put $K_\tau=\bigcup_{i \in I}K_{\tau,i}$ and
\begin{equation*}
 \wh{K}_\tau=\{\xi \in \wh{M} \mid \pi(\xi) \in K_\tau, \|\Pi(\xi)\|\leq 1\}.
\end{equation*}
By $\Delta_\tau$, we denote the Lebesgue number
 of the open cover $(\mc{P}_\tau(\mc{R}_{\tau,i}))_{i \in I_\tau}$
 of $\mc{Z}^-_{\tau}(\mu)$.
Take $\theta_*>0$ and $\nu_*>0$ such that
\begin{equation}
\label{eqn:robust Leb}
 \theta_* \diam \mc{Z}^+_\tau
 +\nu_*^{-1} \diam \mc{Z}^-_\tau(\mu)<\Delta_\tau 
\end{equation}
 for each $\tau=1,2$.
Remark that $f_\tau$ satisfies the $(1,1,1)$-cone condition
 with respect to $(P^h,Q^h)$ on $K_\tau$,
 $\Lambda(K_\tau,f_\tau)$ is a hyperbolic set contained in 
 $\Int K_\tau$,
 and $\wh{f}_\tau$ satisfies
 the $(\theta_*,\mu_*^{-2},\lambda_*)$-cone condition
 with respect to $(\mc{P}_\tau,\mc{Q}_\tau)$ on each $\mc{R}_{\tau,i}$.
Fix an open neighborhood $\mc{V}^-_\tau$ of $E^-_\tau$ such that
\begin{equation*}
 \mc{V}^-_\tau \subset
 \{(z_1,z_2,z_3,Z_{32}) \in \Int \mc{Z}^-_\tau(\mu) \mid
 z_1 \in \Int B^{m_{\tau,1}}(1/4)\}.
\end{equation*}
Let $(f_\sh,\mc{R}_\sh)$ be the pair of the diffeomorphism of $\RR^{2m}$
 and a compact subset of $\wh{M}$
 obtained in Proposition \ref{prop:h-sh}
 for $\theta=\theta_*$, $\nu=\nu_*$,
 $V^-_\tau=\mc{V}^-_\tau$, and $\mc{K}^+_\tau=\mc{Z}^+_\tau$.
\begin{lemma}
\label{lemma:BR intersection}
For any diffeomorphisms $g_1,g_2,g_\sh$ which are $C^2$-close to
 $f_1, f_2, f_\sh$ respectively, 
\begin{equation*}
 \Lambda^s\left(\wh{K}_1,\wh{g_1}\right)
 \cap \wh{g}_\sh\left(\Lambda^s\left(\wh{K}_2,\wh{g}_2\right)
 \right) \neq \emptyset,
\end{equation*}
 where $\wh{g}_1,\wh{g}_2,\wh{g_\sh}$ are the lifts of
 $g_1,g_2,g_\sh$ to $\Gr(\RR^{2m},m)$.
\end{lemma}
\begin{proof}
Put $\lambda_1=\lambda_2=\mu_*^{-2}$ and $\mu_1=\mu_2=\lambda_*$.
Then, $\lambda_1\mu_2=\lambda_2\mu_1=\lambda_*\mu_*^{-2}>1$.
One can check that the blending machines $\mc{B}_1$, $\mc{B}_2$
 for $\wh{f}_1$, $\wh{f}_2$,
 the compact subset $\mc{R}_\sh$ of $\wh{M}$, and
 the map $\wh{f}_\sh$
 satisfy the assumption of Theorem \ref{thm:intersection}
 with $(\theta,\lambda_\tau,\mu_\tau)=(\theta_*,\mu_*^{-2},\lambda_*)$
 and $\lambda_\sh=\mu_\sh=\nu_*$.
\end{proof}

Now, we construct a hyperbolic basic set in $\RR^{2m}$
 which exhibits $C^2$-robust homoclinic tangency of codimension $m$.
We define a diffeomorphism $T:\RR^{2m}\ra \RR^{2m}$ by
 $T(y_1,y_2,y_3,y_4)=(y_4,y_3,y_2,y_1)$ for $y_i \in \RR^{m_{2,i}}$
 $(i=1,2,3,4)$.
By the property of $K_i$ in Proposition \ref{prop:BR},
 the families $(K_{\tau,i})_{i \in I_\tau}$ and 
 $(f_\tau(K_{\tau,i}))_{i \in I_\tau}$ consists of
 mutually disjoint compact subsets in
 $B^{|\vm|-m_{\tau,4}}(4) \times U^+_\tau$
 and $U^-_\tau \times B^{|\vm|-m_{\tau,1}}(4)$ respectively.
Put $K_\sh=B^m(1/4) \times B^m(4)$.
Then, $T(K_\sh)=B^m(4) \times B^m(1/4)$
 and $f_\sh(K_\sh)$ is a subset of $B^m(1/4) \times \RR^m$.
By the choice of $U^\pm_1$ and $U^\pm_2$,
 the compact sets $K_{1,i}$, $T(f_2(K_{2,j}))$, and $T(K_\sh)$
 are mutually disjoint for any $i \in I_1$ and $j \in I_2$.
Similarly, the compact sets $f_1(K_{1,i})$, $T(K_{2,j})$, and $f_\sh(K_\sh)$
 are mutually disjoint for any $i \in I_1$ and $j \in I_2$.
Since all of these compact sets are diffeomorphic to $[-1,1]^{2m}$
 and the maps $f_1$, $f_2$, and $f_\sh$ are orientation preserving,
 there exists a diffeomorphism $f:\RR^{2m} \ra \RR^{2m}$ such that
\begin{equation*}
 f(x)=
\begin{cases}
f_1(x) & x \in K_1,\\
T \circ f_2^{-1} \circ T^{-1}(x) & x \in (T \circ f_2)(K_2)\\
f_\sh \circ T^{-1}(x) & x \in T(K_\sh)
\end{cases}
\end{equation*}
 and the support of $f$ is compact.
Since $\|D(P^h \circ T)v\|=\|DQ^h v\|$
 and  $\|D(Q^h \circ T)v\|=\|DP^h v\|$ for any tangent vector $v$,
 we have
\begin{alignat*}{2}
 DT(C(x,\theta,P^h,Q^h)) & = C(T(x),\theta,Q^h,P^h),\\
 DT(C(x,\theta,Q^h,P^h)) & = C(T(x),\theta,P^h,Q^h), &\hsp
\end{alignat*}
This implies that $T \circ f_2^{-1} \circ T^{-1}$
 satisfies the $(1,1,1)$-cone condition for the splitting $(P^h,Q^h)$
 on $(T \circ f_2)(K_2)$.
Put $K=K_1 \cup (T \circ f_2)(K_2)$.
Then, $f$ satisfies the $(1,1,1)$-cone condition for the splitting $(P^h,Q^h)$
 on $K$.
It is easy to check that
\begin{gather*}
 P^h(K_{1,i})= Q^h(f_1(K_{1,i}))=P^h((T \circ f_2)(K_{2,j}))
 = Q^h(T(K_{2,j}))=B^m(4),\\
 Q^h(K_{1,i}), P^h(f_1(K_{1,i})), Q^h((T \circ f_2)(K_{2,j})),
 P^h(T(K_{2,j})) \subset \Int B^m(4).
\end{gather*}
By Proposition \ref{prop:transitive},
 $\Lambda(K,f)$
 is a topologically transitive hyperbolic set
 in $\Int K$.
Since
\begin{align*}
 \Lambda((f^{-1} \circ T)(K_2),f) 
 & =\bigcap_{n \in \ZZ}f^{-n}((f^{-1} \circ T)(K_2))
 =\bigcap_{n \in \ZZ}f^n(T(K_2))=\Lambda(T(K_2),f^{-1}),
\end{align*}
 the hyperbolic set $\Lambda(K,f)$
 contains $\Lambda(K_1,f) \cup \Lambda(T(K_2),f^{-1})$.
\begin{rmk}
By Remark \ref{rmk:blender K},
 the hyperbolic set $\Lambda(K_1,f)=\Lambda(K_1,f_1)$ is a $cu$-blender for $f$.
Similarly, $\Lambda(T(K_2),f)=T(\Lambda(K_2,f_2^{-1}))$
 is a $cu$-blender of $f^{-1}=T \circ f_2 \circ T^{-1}$,
 and hence, it is a $cs$-blender for $f$.
Since $\Lambda(K,f)$ is a topologically transitive
 hyperbolic set which contains both $\Lambda(K_1,f)$ and $\Lambda(T(K_2,f))$,
 it is a {\it double blender} for $f$ {\it i.e.},
 a $cu$- and $cs$-blender simultaneously.
Such a blender was introduced in \cite{NP12}
 and used in \cite{ACW-pre,BR17,BR18}.
\end{rmk}
\begin{figure}
\begin{center}
 \includegraphics[scale=0.6]{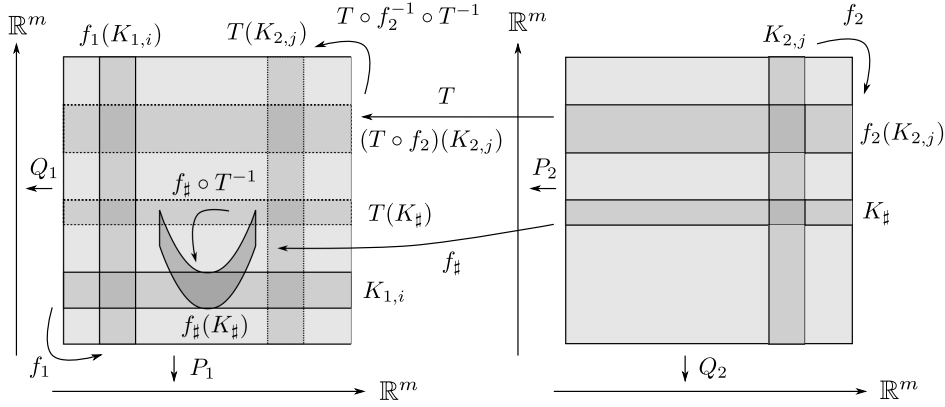} 
\end{center}
 \caption{Proof of Theorem B}
\label{fig:Cantor}
\end{figure}

Let $\wh{f}$ and $\wh{T}$ be the lifts of $f$ and $T$ to $\Gr(\RR^{2m},m)$.
Then, 
\begin{equation*}
\wh{f}(\xi)=
\begin{cases}
 \wh{f}_1(\xi) & (\xi \in \wh{K}_1),\\
 \wh{T} \circ \wh{f}_2^{-1} \circ \wh{T}^{-1}(\xi)
 & (\xi \in \wh{T}(\wh{f}_2(\wh{K}_2)).
\end{cases}
\end{equation*}
Put $\wh{K}_\sh=\{\xi \in \wh{M} \mid \pi(\xi) \in K_\sh,
 \|\Pi(\xi)\|\leq 1\}$.
Then,
\begin{equation*}
 \mc{R}_\sh \subset (\mc{P},\mc{Q})^{-1}(\mc{V}^-_2 \times \mc{Z}^+_2)
 \subset \wh{K}_\sh,
\end{equation*}
 and hence, $\wh{f}= \wh{f}_\sh \circ \wh{T}^{-1}$ on $\wh{T}(\mc{R}_\sh)$.
By the persistence of a locally maximal hyperbolic set,
 there exist a $C^1$-neighborhood $\cU_0$ of $f$ such that
 $\Lambda(K,g)$ is a topologically transitive hyperbolic set of $g$
 and $\Int K$ contains $\Lambda(K_1,g)$ and $\Lambda(T(K_2),g^{-1})$
 for any $g \in \cU_0$.
For any diffeomorphism $g$ sufficiently $C^2$-close $f$,
 the maps $\wh{g}$,
 $\wh{T}^{-1} \circ \wh{g}^{-1} \circ \wh{T}$,
 and $\wh{g} \circ \wh{T}$ are $C^1$-close to
 $\wh{f}_1$, $\wh{f}_2$, and $\wh{f}_\sh$
 on $\wh{K}_1$, $\wh{K_2}$, and $\wh{K}_\sh$ respectively.
By Lemma \ref{lemma:BR intersection},
 there exists a $C^2$-neighborhood $\cU$ of $f$
 such that $\cU \subset \cU_0$ and
\begin{equation}
\label{eqn:robust intersection}
 \Lambda^s(\wh{K}_1,\wh{g}) \cap
 (\wh{g} \circ \wh{T})
 (\Lambda^s(\wh{K}_2,\wh{T}^{-1} \circ \wh{g}^{-1} \circ \wh{T}))
 \neq \emptyset
\end{equation}
 for any $g \in \cU$,
 where $\wh{g}$ is the lift of $g$ to $\Gr(\RR^{2m},m)$.
Since
\begin{align*}
\wh{g}_\sh(\Lambda^s(\wh{K}_2,\wh{T}^{-1} \circ \wh{g}^{-1} \circ \wh{T}))
 & =\wh{g} \circ \wh{T}\left(\bigcap_{n \geq 0}
  (\wh{T}^{-1} \circ\wh{g}^{-1} \circ \wh{T})^{-n}(\wh{K}_2)\right)\\
 & = \wh{g}\left(\bigcap_{n \geq 0}\wh{g}^n(\wh{T}(\wh{K}_2))\right)
   =\wh{g}(\Lambda^s(\wh{T}(\wh{K}_2),\wh{g}^{-1})),
\end{align*}
 this implies that
\begin{equation}
\label{eqn:lift intersection}
 \Lambda^s(\wh{K}_1,\wh{g}) \cap
 \wh{g}(\Lambda^s(\wh{T}(\wh{K}_2),\wh{g}^{-1})
 \neq \emptyset
\end{equation}
 for any $g \in \cU$.
By Proposition \ref{prop:lift W-s}, we have
\begin{align*}
 \Lambda^s(\wh{K}_1,\wh{g})
 & \subset \{T_q W^s(p,g) \mid p \in \Lambda(K_1,g), q \in W^s(p,g_1)\},\\
\wh{g}(\Lambda^s(\wh{T}(\wh{K}_2),\wh{g}^{-1}))
 & \subset \{T_q W^s(p,g^{-1}) \mid p \in \Lambda(T(K_2),g^{-1}),
 q \in W^s(p,g^{-1})\}\\
 & = \{T_q W^u(p,g) \mid
  p \in \Lambda(T(K_2),g^{-1}), q \in W^u(p,g)\}.
\end{align*}
By (\ref{eqn:robust intersection}),
 for any given $g \in \cU$
 there exist $p_1 \in \Lambda(K_1,g)$,
 $p_2 \in \Lambda(T(K_2),g^{-1})$,
 and $q \in W^s(p_1,g) \cap W^u(p_2,g)$
 such that $T_q W^s(p_1,g)=T_q W^u(p_2,g)$.
Since $\Lambda(K_1,g)$ and $\Lambda(T(K_2),g^{-1})$ are $g$-invariant subset
 in $\Int K$, they are contained in  $\Lambda(K,g)$.
This implies that
 the topologically transitive hyperbolic set $\Lambda(K,g)$
 exhibits homoclinic tangency of codimension $m$ at $q$.
Therefore, $f$ exhibits $C^2$-robust homoclinic tangency
 of codimension $m$.
This finish the proof of Theorem \ref{thm:ThmB'},
 and hence, of Theorem B.
\begin{rmk}
By Remark \ref{rmk:blender Gr},
 $\Lambda(\wh{K}_1,\wh{g})$ and $\Lambda(\wh{T}(\wh{K}_2),\wh{g})$
 are hyperbolic invariant sets of $\wh{g}$ if $g$ is sufficiently $C^2$-close
 to $f$.
The condition (\ref{eqn:lift intersection}) implies
 that there exists a heteroclinic point between these hyperbolic sets.
\end{rmk}

%
%


\begin{thebibliography}{99}
\bibitem{Ar02}
 Z.Arai,
 Tangencies and the Conley index.
 Ergodic Theory Dynam. Systems 22 (2002), no. 4, 973--999.

\bibitem{AM-pre}
H.Ara\'VJ and C.G.Moreira,
 Stable intersections of conformal Cantor sets.
 \url{arXiv:1910.03715}.

\bibitem{ACW-pre}
 A. Avila, S.Crovisier, and A.Wilkinson,
 $C^1$ density of stable ergodicity.
 \url{arXiv:1709.04983}.

\bibitem{Ba}
L.Barreira,
 Dimension and recurrence in hyperbolic dynamics.
 Progress in Mathematics, 272. Birkh\"{o}user Verlag, Basel, 2008.
 xiv+300 pp.

\bibitem{BFMS2017}
P.G.Barrientos, A.Fakhari,D.Malicet, and A.Sarizadeh,
 Expanding actions: minimality and ergodicity.
 Stoch. Dyn. 17 (2017), no.4, 1750031, 
 \url{DOI: 10.1142/S0219493717500319}.

\bibitem{BKR14}
P.G.Barrientos, Y.Ki, and A.Raibekas,
 Symbolic blender-horseshoes and applications. 
 Nonlinearity 27 (2014), no. 12, 2805--2839.

\bibitem{BR17}
P.G.Barrientos and A.Raibekas,
 Robust tangencies of large codimension,
 Nonlinearity {\bf 30} (2017), 4369--4409.

\bibitem{BR18}
 P.G.Barrientos and A.Raibekas,
 Robustly non-hyperbolic transitive symplectic dynamics.
 Discrete Contin. Dyn. Syst. 38 (2018), no. 12, 5993--6013.

\bibitem{BR20}
 P.G.Barrientos and A.Raibekas,
 Robust degenerate unfoldings of cycles and tangencies.
 J. J Dyn Diff Equat (2020). \url{DOI:10.1007/s10884-020-09857-0}.

\bibitem{BP-pre}
P.G.Barrientos and S.A.P\'erez,
 Robust heteroclinic tangencies, \url{arXiv:1903.00264}.

\bibitem{Be16}
P.Berger, Generic family with robustly infinitely many sinks.
 Invent. Math. 205 (2016), no. 1, 121--172.
 (Correction: Invent. Math. 218 (2019), no. 2, 649--656)

\bibitem{Bi-pre}
S.Biebler,
 A complex gap lemma, \url{arXiv:1810.02544}.

\bibitem{BBD16}
J.Bochi, C.Bonatti, and L.J.L\'iaz, 
 Robust criterion for the existence of nonhyperbolic ergodic measures.
 Comm. Math. Phys. 344 (2016), no. 3, 751--795. 

\bibitem{BD96}
C.Bonatti and L.J.D\'iaz,
 Persistent nonhyperbolic transitive diffeomorphisms.
 Ann. of Math. (2) 143 (1996), no. 2, 357--396. 

\bibitem{BD08}
C.Bonatti and L.J.D\'iaz,
 Robust heterodimensional cycles and $C^1$-generic dynamics.
 J. Inst. Math. Jussieu 7 (2008), no. 3, 469--525. 

\bibitem{BD12}
C.Bonatti and L.J.D\'iaz,
 Abundance of $C^1$-robust homoclinic tangencies.
 Trans. Amer. Math. Soc. 364 (2012), no. 10, 5111--5148.

\bibitem{BDK12}
C.Bonatti, L.J.D\'iaz, and S.Kiriki,
Stabilization of heterodimensional cycles.
Nonlinearity 25 (2012), no. 4, 931--960. 

\bibitem{BDV}
C.Bonatti, L.J.D\'iaz, and M.Viana,
Dynamics beyond uniform hyperbolicity.
 A global geometric and probabilistic perspective.
 Encyclopaedia of Mathematical Sciences, 102.
 Springer-Verlag, Berlin, 2005. xviii+384 pp.

\bibitem{Bu93}
G.T.Buzzard, Stably interesting Julia sets of polynomials.
 C. R. Acad. Sci. Paris S\'er. I Math. 317 (1993), no. 11, 1013--1018.

\bibitem{Bu97}
G.T.Buzzard,
 Infinitely many periodic attractors for holomorphic maps of 2 variables.
Ann. of Math. (2) 145 (1997), no. 2, 389--417. 

\bibitem{BCF18}
J.Buzzi, S.Crovisier, and T.Fisher,
 The entropy of $C^1$-diffeomorphisms without a dominated splitting. 
 Trans. Amer. Math. Soc. 370 (2018), no. 9, 6685--6734. 

\bibitem{Ca19}
T.Catalan,
 A link between topological entropy and Lyapunov exponents.
 Ergodic Theory Dynam. Systems 39 (2019), no. 3, 620--637.

\bibitem{DGS15}
D.Damanik, A.Gorodetski, and B.Solomyak,
 Absolutely continuous convolutions of singular measures
 and an application to the square Fibonacci Hamiltonian.
Duke Math. J. 164 (2015), no. 8, 1603--1640. 


\bibitem{GTS93-1}
 S.V.Gonchenko,  D.V.Turaev, and L.P.Shil'nikov,
 On the existence of Newhouse regions in a neighborhood
 of systems with a structurally unstable homoclinic Poincar\'e curve
 (the multidimensional case).
 (Russian) Dokl. Akad. Nauk 329 (1993), no. 4, 404--407;
 translation in Russian Acad. Sci. Dokl. Math. 47 (1993), no. 2, 268--273.

\bibitem{GTS93-2}
 S.V.Gonchenko,  D.V.Turaev, and L.P.Shil'nikov,
 Dynamical phenomena in multidimensional systems
 with a structurally unstable homoclinic Poincar\'e curve.
 (Russian) Dokl. Akad. Nauk 330 (1993), no. 2, 144--147;
 translation in Russian Acad. Sci. Dokl. Math. 47 (1993), no. 3, 410--415

\bibitem{KH}
 A.~Katok and B.~Hasselblatt,
 Introduction to the modern theory of dynamical systems.
 Encyclopedia of Mathematics and its Applications,
 54. Cambridge University Press, Cambridge, 1995.

\bibitem{KS12}
S.Kiriki and T.Soma,
 $C^2$-robust heterodimensional tangencies. 
Nonlinearity 25 (2012), no. 12, 3277--3299.

\bibitem{Mo96}
C.G.Moreira,
 Stable intersections of Cantor sets and homoclinic bifurcations.
 Ann. Inst. H. Poincar\'e Anal. Non Lin\'eaire 13 (1996), no. 6, 741--781. 

\bibitem{Mo98}
C.G.Moreira,
Sums of regular Cantor sets, dynamics and applications to number theory.
International Conference on Dimension and Dynamics (Miskolc, 1998).
Period. Math. Hungar. 37 (1998), no. 1--3, 55--63. 

\bibitem{Mo11}
C.G.Moreira,
 There are no $C^1$-stable intersections of regular Cantor sets.
 Acta Math. 206 (2011), no. 2, 311--323. 

\bibitem{MY01}
C.G.T.de A.Moreira and J.-C.Yoccoz,
 Stable intersections of regular Cantor sets with large Hausdorff dimensions.
Ann. of Math. (2) 154 (2001), no. 1, 45--96. 

\bibitem{NP12}
M.Nassiri and E.R.Pujals,
 Robust transitivity in Hamiltonian dynamics.
Ann. Sci. \'Ec. Norm. Sup\'er. (4) 45 (2012), no. 2, 191--239. 

\bibitem{Ne74}
S.E.Newhouse, 
Diffeomorphisms with infinitely many sinks. Topology {\bf 13} (1974), 9--18. 

\bibitem{Ne79}
S.E.Newhouse, 
The abundance of wild hyperbolic sets and nonsmooth stable sets
 for diffeomorphisms.
 Inst. Hautes \'Etudes Sci. Publ. Math. No. 50 (1979), 101--151.

\bibitem{PT87}
J. Palis and F. Takens,
Hyperbolicity and the creation of homoclinic orbits.
Ann. of Math. (2) 125 (1987), no. 2, 337--374. 

\bibitem{PT93}
J. Palis and F. Takens,
 Hyperbolicity and sensitive-chaotic dynamics at homoclinic bifurcations.
 Cambridge University Press, 1993.

\bibitem{PV94}
J.Palis and M.Viana,
 High dimension diffeomorphisms displaying
 infinitely many periodic attractors.
 Ann. of Math. (2) 140 (1994), no. 1, 207--250. 

\bibitem{PY94}
J.Palis and J.-C.Yoccoz,
 Homoclinic tangencies for hyperbolic sets of large Hausdorff dimension.
 Acta Math. 172 (1994), no. 1, 91--136. 

\bibitem{Sh}
 M. Shub,
Global stability of dynamical systems, Springer-Verlag, 1987.

\bibitem{Ur95}
R.Ures,
Abundance of hyperbolicity in the $C^1$ topology.
Ann. Sci. \'Ecole Norm. Sup. (4) 28 (1995), no. 6, 747--760. 


\end{thebibliography}
\end{document}